\documentclass[titlepage,a4paper]{article}
\usepackage[T1]{fontenc}
\usepackage[latin1]{inputenc}
\usepackage{amsmath}
\usepackage{amssymb}

\usepackage{amsthm}

\newtheorem{ex}{Example}
\newtheorem{Def}{Definition}

\newtheorem{rem}{Remark}
\newtheorem{Th}{Theorem}
\newtheorem{res}{Result}
\newtheorem{lem}{Lemma}

\newtheorem{cor}{Corollary}

\newcommand{\PP}{\mathbb{P}}
\newcommand{\EE}{\mathbb{E}}

\newcommand{\RR}{\mathbb{R}}

\newcommand{\NN}{\mathbb{N}}
\newcommand{\ZZ}{\mathbb{Z}}

\newcommand{\cG}{{\cal G}}

\newcommand{{{\cadlag}}}{c\`adl\`ag}

\newcommand{\esssup}{\mathrm{ess \ sup} \ }
\newcommand{\id}{\mathrm{id}}

\begin{document}

\sffamily
\sloppy
\pagestyle{plain}

\renewcommand{\baselinestretch}{1}

\title{Error bounds and convergence for American put option pricing based on translation-invariant Markov chains}
\author{Frederik S Herzberg \thanks{Abteilung f\"ur Stochastik, Institut f\"ur Angewandte Mathematik, Universit\"at Bonn, D-53115 Bonn, Germany ({\tt herzberg@wiener.iam.uni-bonn.de})} \thanks{Mathematical Institute, University of Oxford, Oxford OX1 3LB, England} }
\date{}

\maketitle

\renewcommand{\baselinestretch}{1}

\begin{abstract} Consider a discrete finite-dimensional, Markovian market model. In this setting, discretely sampled American options can be priced using the so-called ``non-recombining'' tree algorithm. By successively increasing the number of exercise times, the American option price itself can be computed; for combinatorial reasons, we shall consider a recursive algorithm that doubles the number of exercise times at each recursion step. First we prove, by elementary arguments, error bounds for the first order differences in this recursive algorithm. From this, bounds on the higher order differences can be obtained using combinatorial arguments that are motivated by the theory of rough paths. We shall obtain an explicit $L^1(C)$ convergence estimate for the recursive algorithm that prices a discretely sampled American $\max$-put option (on a basket of size $d$) at each recursion step, $C$ belonging to a certain class of compact subset of $\RR^d$, in under the assumption of sufficiently small volatilities. In case $d=1$, $L^1(C)$-bounds for an even more natural choice of $C$ will be derived.
\end{abstract}

\noindent

\label{cubature_dyadic}

\section{Introduction}

An American option is a derivative security on a number of assets which can be exercised at any time before maturity. In case the maturity is a finite number $T$, which we shall suppose throughout this paper, this option will be called ``non-perpetual''. If the set of possible exercise times is discrete, these options are also referred to as {\em discretely sampled American} (alternatively shorter: {\em Bermudan}) options. A holder of such an options will only exercise this option if the expected gain she can draw from holding on to the option falls short of what she would get from exercising the option in the current situation of the market. In this paper, we shall look at a pricing algorithm that is based on a market model that is sometimes referred to as a ``non-recombining tree'' (eg Glasserman \cite{Gla}) or simply ``lattice'' (eg Kargin \cite{Karg}). This numerical procedure differs from a Markov chain-Monte Carlo algorithm in that it is intrinsically deterministic (cf again Glasserman \cite{Gla}), and has already been studied by Heath, Jarrow and Morton \cite{HJM90}. The nodes and weights of this market model might, for example, be derived from a cubature formula for a Gaussian measure. This could be justified by regarding the discrete market model merely as an approximation of a time-continuous Black-Scholes model rather than taking the discrete market model as a model in its own right -- which may be difficult to vindicate, given that for completeness the number of branches that leave from each node must equal the number of assets in the basket plus one (the bond).

\subsection{Motivating Bermudan option pricing based on cubature}

If the translation-invariant Markov chain of the market model is indeed derived from a cubature formula ``with few points'' \cite{V} for a symmetric measure, its set of increments will also be, whilst being ``asymmetric'' \cite{V}, highly regular. Therefore the -- from a numerical perspective undesirable -- attribute ``non-recombining'' in the term ``non-recombining tree'' will not be applicable any longer. To be more precise, one can note

\begin{rem} \label{treeispolynomial} Suppose $\left\{\xi_1,\dots,\xi_m\right\}\subset\RR^d$ and $\left\{\alpha_1,\dots,\alpha_m\right\}$ are the sets of cubature points and corresponding weights, respectively, from the cubature formula for the integration of degree $5$ polynomials with respect to a standard Gaussian measure of dimension $d=3k-2$ (for some arbitrary $k\in\NN$) from Victoir's example \cite[5.1.1]{V}. Consider $r\geq 0$ and $\sigma,\mu,T>0$. Then the pointwise recursion defined by \begin{eqnarray}\nonumber \tilde V_T&=&g,\\ \label{recursion} \nonumber \forall k\in\NN\cap\left[0,\frac{T}{h}\right]&& \\ \tilde V_{(k-1)h}(\cdot)&=&\max\left\{e^{-rh}\sum_{j=1}^m\alpha_j \tilde V_{kh}(\cdot +\mu h+h^{1/2}\sigma\cdot\xi_j), g(\cdot)\right\} \end{eqnarray} (provided $T$ is an integer multiple of $h>0$) will only have to apply the function $g$ to a number of points in $\RR^d$ that grows polynomially in $\frac{1}{h}$ for $h>0$.
\end{rem}
\begin{proof} The cubature points of the cubature formulae referred to in the Theorem form a finite subset of $\sqrt[4]{3}\{0,\pm 1\}^d$. Sums of length $\frac{T}{h}$ (provided this fraction is an integer) of the cubature points are therefore always elements of $\sqrt[4]{3}\left(\ZZ^d\cap\left\{\left|\cdot\right|\leq \frac{1}{h}\right\}\right)$ (and this set has only $\left(2\frac{1}{h}\right)^d$ elements), and the points used in the recusion formula stated above are comprised of a subset of $h^{1/2}\cdot\sqrt[4]{3}\left(\ZZ^d\cap\left\{\left|\cdot\right|\leq \frac{1}{h}\right\}\right)+ h \cdot\mu \left\{0,\dots,\frac{1}{h}\right\}+\xi_0$.
\end{proof}

However, this is not the only recombination that can be accomplished in the case where $d=3k-2$ :

\begin{rem} Let, for the sake of simplicity, $k=3$ and thus $d=7$, and consider the cubature formula of degree $5$ for the Gaussian measure found by Victoir \cite[Example 5.1.1]{V}, the cubature points being $\{x_0\}\cup\sqrt[4]{3}\cdot\cG_3X_1\subseteq\RR^7$ (in Victoir's notation), where $x_0=0\in\RR^7$ and $X_1$ is some seven-element subset of $\{0,1\}^7$ on which the group $\cG_3$, which is the group generated by the three reflections with respect to the fourth, sixth, and seventh coordinate axis (ie the group generated by the maps $(x_1,\dots,x_7)\mapsto (x_1,\dots,x_3,-x_4,x_5,\dots,x_7)$, $(x_1,\dots,x_7)\mapsto (x_1,\dots,x_5,-x_6,x_7)$, and $(x_1,\dots,x_7)\mapsto (x_1,\dots,x_6,-x_7)$), acts. 

The recursion of Remark \ref{treeispolynomial}, can be regarded formally as a tree where at each node exactly one branch leaves for each element of the set $\{0\}\cup\sqrt[4]{3}\cdot\cG_3X_1$, $0\in\RR^d$ being the root of the tree

In order to find and eliminate those branches of the tree that are computed ``wastefully'', one can divide the sums (of length $\frac{1}{h}$) of the cubature points by $\sqrt[4]{3}$ (whence one obtains a subset of $\ZZ^7$) and consider them coordinate-wise modulo $2$. Then one is dealing with elements of the vector space $\left(\ZZ/2\ZZ\right)^d$. The coordinate-wise projection of the $\frac{1}{\sqrt[4]{3}}$-multiple of the set of cubature points $\{x_0\}\cup \sqrt[4]{3}\cdot\cG_3X_1$ into the vector space $\left(\ZZ/2\ZZ\right)^7$ now contains only eight points (instead of $57$ as before). 

Thus, using basic linear algebra in a $7$-dimensional $\ZZ /2\ZZ$-vector space, we are easily able to classify the non-trivial zero representations from elements of the projected cubature points. 

Perceiving $X_1$ as a $7$-element subset of $\left(\ZZ /2\ZZ\right)^7$, we see that $(x)_{x\in X_1}$ is an invertible $\left(\ZZ /2\ZZ\right)^{7\times 7}$-matrix. Therefore we cannot expect any recombination from representations of zero by nontrivial linear combinations of elements of $X_1\subset \left(\ZZ /2\ZZ\right)^7$. Moreover, the fact that $A:=(x)_{x\in X_1}$ is invertible, shows that $x_0 ={ 0}$ can only be written trivially as a sum of elements of $X_1$. Hence we have shown that we exploit symmetries optimally if we use: (i) the commutativity of $(\RR^d, +)$; (ii) the obvious symmetries due to the construction of the cubature formulae by means of the action of a reflection group on certain points; (iii) the fact that addition of $x_0$ does nothing at all.
\end{rem}

Hence, in terms of complexity reduction, cubature based on the cubature points found by Victoir is clearly promising.

\subsection{Problem formulation and notation}

In order to introduce the underlying discrete market model, let us consider an arbitrary but fixed translation-invariant finite-state Markov chain $P:=\left(P_t\right)_{t\in I}$ with state space $\RR^d$ (for $d\in\NN$) where $I= h\NN_0$ for some real number $h>0$, as well as a real number $T\in h\NN_{>1}\subset I$ (the time horizon, or maturity), a real number $r>0$ (the discount rate), a continuous function $\bar f:\RR^d\rightarrow [0,+\infty)$ that is monotone in each coordinate (defining the contingent claim as a function of the {\em logarithmic} prices of the assets in the basket), a nonnegative real number $K\geq 0$ (the strike price), and let $g:\RR^d\rightarrow\RR$ be measurable. (Often we shall assume $$g:=K- \bar f,$$ such that $g\vee 0$ is the payoff function for the corresponding put.) We will also define a family of maps $B_t:L^0\left(\RR^d,[0,+\infty)\right)\rightarrow L^0\left(\RR^d,[0,+\infty)\right)$, $t\in I$, by $$\forall t\in I \quad B_t:f\mapsto \max\left\{e^{-rt}P_tf,g\right\}= \left(e^{-rt}P_tf\right)\vee g.$$ (Note that $B_tf$ will always be nonnegative for $f\geq 0$ -- hence, for all $f\geq 0$, $B_tf\geq g\vee 0$.) Furthermore, we shall denote by $\left\{y_1^{(t)},\dots, y_{m^{\frac{t}{h}}}^{(t)}\right\}$ the set of increments of the translation-invariant Markov chain ($y_1^{(t)},\dots, y_{m^{\frac{t}{h}}}^{(t)}$ not necessarily mutually distinct) -- and by $\left\{\alpha^{(t)}_1,\dots,\alpha^{(t)}_{m^\frac{t}{h}}\right\}\subset(0,1]$ the set of the corresponding transition probabilities, implying in particular $\sum_{i=1}^{m^{\frac{t}{h}}}\alpha_{i}^{(t)}=1$ for all $t\in I$. Hence $$\forall t\in I\quad P_t:f\mapsto \sum_{i=1}^{m^{\frac{t}{h}}}\alpha_{i}^{(t)}f\left(\cdot+y_i^{(t)}\right)$$ and \begin{equation}\left\{y_1^{(t)},\dots, y_{m^{\frac{t}{h}}}^{(t)}\right\}=\left\{\sum_{\ell=1}^{t/h} y_{k_\ell}^{(h)} \ : \ k_1,\dots,k_{\ell}\in \{1,\dots, m\}\right\}.\label{chapman}\end{equation}

If there are $N$ equidistant exercise times before maturity, then the expected payoff of the corresponding Bermudan option with payoff function $g\vee 0$ for $g=K-f$ would be $\left(B_{T/N}\right)^{\circ N}(g\vee 0)$, the motivation being that at each time-step the holder option will decide whether the expected gain from holding on to the option exceeds what she would get from exercising the option now. For an approximation of the American (rather than the Bermudan) price, one approach would be to choose a dyadic partition of $[0,T]$ by choosing powers of $2$ as our $N$ in the preceding expression. This makes the sequence of functions $\left(\left(B_{T\cdot 2^{-n}}\right)^{\circ 2^n}(g\vee 0)\right)_{n\in\NN_0}$ pointwise monotonely increasing, and our aim in this paper is to establish error bounds and/or convergence estimates for this sequence.

Before bringing the introduction to a close, let us set up notation. The Lebesgue measure on $\RR$ shall be denoted by $\lambda$, $\lambda^d$ will be the Lebesgue measure on $\RR^d$. The operators $\max$ and $\min$ when applied to subsets of $\RR^d$ will be understood to be taken componentwise. Analogously, we will interpret the relations $\leq$ and $\geq $ componentwise on $\RR^d$. Thus, eg the assertion $\left(\max_{k\in\{1,\dots,m\}} \left(z_k\right)_j\right)_{j\in\{1,\dots,d\}}\geq 0$ componentwise (for $z_1,\dots,z_m\in\RR^d$) will be written just $\max_{k\in\{1,\dots,m\}} z_k \geq 0$.

For convenience, we allow all $L^p$-norms (including the $L^\infty$ norm) of measurable functions to take values in the interval $[0,+\infty]$, thereby extending the domain for each of the $L^p$-norm to $L^0$, the vector lattice of measurable functions. Furthermore, any functions occurring in this paper will be assumed to be measurable. Thus, eg the relation $f_0\geq f_1$ should be read as shorthand for $f_0\in L^0\left(\lambda^d\right)\cap\left\{\cdot \geq f_1\right\}$ for all functions $f_0,f_1$; analogously for the relation $f_0\leq f_1$.

We will use the operation $\vee$ in such a way that it is applied prior to $+$, but only after $P_s$ and multiplication with other functions or constants have taken place: $$C\cdot P_sf_0\vee f_3\cdot f_1 +f_2 =\max\left\{C\cdot P_sf_0,f_3\cdot f_1\right\}+f_2.$$

\subsection{Main results}

The main results of this paper are the following two $L^1$ estimates:

\begin{res}[Theorem \ref{cubatureBermudanconv}] Suppose $d=1$ and $\bar f=\exp$, as well as $g=K-\bar f$. Under these assumptions there is a $\gamma_1$ such that $P_t\bar f={\gamma_1}^t\bar f$ for every $t\in I$, and let us suppose this $\gamma_1\in(0,e^r]$. Assume furthermore $y_i^{\left(h\right)}\geq 0$ componentwise for all ${i\in\{1,\dots,m\}}$. Then there exists a real number $D>0$ such that for all $N>M\in\NN$, $s\in(0,T]\cap \left(2^N\cdot I\right)$ and monotonely decreasing $f\geq g\vee 0$, one has \begin{eqnarray*}&& \left\|\left(B_{s\cdot 2^{-N}}\right)^{\circ\left(2^{N}\right)}f - \left(B_{s\cdot 2^{-M}}\right)^{\circ\left(2^M\right)}f\right\|_{L^1\left(\left\{P_h(g\vee 0)>P_hg\right\}\right)}\\ &\leq& D\cdot {s}^2\cdot{2}^{-M}\left(1- 2^{-(N-M-1)}\right) \\ &\leq& D\cdot {s}^2\cdot{2}^{-M}.\end{eqnarray*} 
$D$ is computed explicitly in the statements of Lemma \ref{estimateonE} and Theorem \ref{estimateBs/2-BsonE}.
\end{res}

\begin{res}[Theorem \ref{cubatureBermudanconv_barf=min}] Suppose $\bar f=\max_{j\in\{1,\dots,d\}}\exp\left((\cdot)_j\right)$ and consider a compact set $B$ such that $B-(\ln K)_{j=1}^d\subseteq\left[-{R},R\right]^d$. Assume that $y_{i}^{(h)}\geq 0$ for all $i\in\{1,\dots,m\}$ and $0\in\left\{y_i^{(h)} \ : \ i\in\{1,\dots,m\}\right\}$. If one now defines $$D=\left(\left(\ln\gamma_1-r\right)\tilde{D} +rK+ C_0  \right)\cdot R^{d-1}\cdot \frac{\sum_{j=1}^d\max_{i} \left(\left(y_i^{(h)}\right)_j\vee 0\right)}{h} $$ (with $C_0$ and $\tilde D$ as defined in Lemma \ref{estimateonE} and the $\gamma_1$ of Lemma \ref{barf=max_gamma1}), then one has for all $N>M\in\NN$, $s\in(0,T)\cap \left(2^N\cdot I\right)$ and $f\geq g\vee 0$,
\begin{eqnarray*}&&\left\|\left(B_{s\cdot 2^{-N}}\right)^{\circ\left(2^{N}\right)}f - \left(B_{s\cdot 2^{-M}}\right)^{\circ\left(2^M\right)}f\right\|_{L^1\left(\left\{P_h(g\vee 0)>P_hg\right\} \cap \bigcap_{\ell}\left(B-y_\ell^{\left(s\right)}\right)\right)} \\ &\leq& D\cdot {s}^2\cdot{2}^{-M}\left(1- 2^{-(N-M-1)}\right) \\ &\leq& D\cdot {s}^2\cdot{2}^{-M}\longrightarrow 0 \text{ as }M\rightarrow\infty. \end{eqnarray*}
\end{res}

\section{Monotonicity of the algorithm, bounds on first order differences, and the family $E^\cdot$}

First of all, we state the pointwise montonocity of approximate American option pricing based on dyadic partitions:

\begin{lem} The sequence $\left(\left(B_{T\cdot 2^{-n}}\right)^{\circ 2^n}f\right)_{n\in\NN_0\cap\left\{T\cdot 2^{-\cdot}\in I\right\} }$ is pointwise monotonely increasing for all functions $f:\RR^d\rightarrow \RR$. Furthermore, if there exists a function $\tilde g\geq g\vee 0$ such that $\tilde g$ is {\em $e^{-r\cdot}P_\cdot$-harmonic} (ie $e^{-rh}P_h\tilde g=\tilde g$) and $f\leq \tilde g$, then for all $n\in\NN_0$, $\left(B_{T\cdot 2^{-n}}\right)^{\circ 2^n}f\leq \tilde g$.
\end{lem}
\begin{proof} Consider $f:\RR^d\rightarrow \RR$ and $n\in\NN_0$ such that $T\cdot 2^{-(n+1)}\in I=h\NN_0$. Then $$\left(B_{T\cdot 2^{-(n+1)}}\right)^{\circ 2^{n+1}}f=\left(\left(B_{T\cdot 2^{-n+1}}\right)^{\circ 2}\right)^{\circ 2^n}f$$ and by the monotonicity of the operators $P_s$ for $s\in I$, \begin{eqnarray*}\left(B_{T\cdot 2^{-n+1}}\right)^{\circ 2}&=&e^{-rT \cdot 2^{-(n+1)}}P_{T \cdot 2^{-(n+1)}}\left(e^{-rT \cdot 2^{-(n+1)}}P_{T \cdot 2^{-(n+1)}}(\cdot)\vee g\right)  \vee g\\ &\geq &e^{-rT \cdot 2^{-(n+1)}}P_{T \cdot 2^{-(n+1)}}\left(e^{-rT \cdot 2^{-(n+1)}}P_{T \cdot 2^{-(n+1)}}(\cdot)\right)\vee g\\ &=& e^{-rT \cdot 2^{-n}}P_{T \cdot 2^{-n}}(\cdot)\vee g=B_{T\cdot 2^{-n}},\end{eqnarray*} where the last line is a consequence of the Chapman-Kolmogorov equation. This completes the proof for the monotonicity of the sequence $\left(\left(B_{T\cdot 2^{-n}}\right)^{\circ 2^n}f\right)_{n\in\NN_0\cap\left\{T\cdot 2^{-\cdot}\in I\right\}}$. 

Now suppose there exists such a function $\tilde g$ as in the statement of the Lemma. Then $e^{-rs}P_s\tilde g=\tilde g$ for all $s\in I$ and therefore $B_s\tilde g=\tilde g$ for all $s\in I$. Also, the map $B_s$ is monotone in the sense that $g_0\leq g_1$ always implies $B_sg_0\leq B_sg_1$ (because it is the composition of two monotone maps: $e^{-rs}P_s$ and $\cdot\vee g$) for all $s\in I$. Thus we see that for all $f\leq\tilde g,$ $$\left(B_{T\cdot 2^{-n}}\right)^{\circ 2^n}f\leq \left(B_{T\cdot 2^{-n}}\right)^{\circ 2^n}\tilde g =\tilde g.$$
\end{proof}

\begin{lem} \label{B_tcontractsifgeq g} For all measurable functions $f_1\geq f_0\geq g\vee 0$, as well as for all $t\in I$ and $p\in\{1,\infty\}$ one has \begin{eqnarray*} \left\|B_tf_1-B_tf_0\right\|_{L^p\left(\lambda^d\right)} &\leq& e^{-rt} \left\|f_1-f_0\right\|_{L^p\left(\lambda^d\left[\left\{e^{-rt}P_tf_1>g\right\}\cap\cdot\right]\right)} \\ &\leq& e^{-rt} \left\|f_1-f_0\right\|_{L^p\left(\lambda^d\right)}  \end{eqnarray*} (with the usual convention that $x\leq +\infty$ for all $x\in\RR\cup\{\pm\infty\}$).
\end{lem}
\begin{proof} The map $B_t$ is monotone. Thus we have \begin{eqnarray*} \left\{B_tf_1=g\right\}&=& \left\{B_tf_0\leq B_tf_1=g\right\} \\ &=&  \left\{g\vee 0 \leq B_tf_0\leq B_tf_1=g\right\} = \left\{B_tf_0=g\right\}\cap  \left\{B_tf_1=g\right\}\\ &\subseteq&  \left\{B_tf_1-B_tf_0=0\right\} \end{eqnarray*} for $f_1\geq f_0\geq g\vee 0$.
Since $B_tf_1\geq g$, this implies \begin{eqnarray*}0\leq B_tf_1-B_tf_0&=& \chi_{\left\{B_tf_1>g\right\}}\left(e^{-rt}P_tf_1\vee g - e^{-rt}P_tf_0\vee g\right) \\ &=& \chi_{\left\{B_tf_1>g\right\}}\left(e^{-rt}P_tf_1-e^{-rt}P_tf_0\vee g\right)\\ &\leq & \chi_{\left\{e^{-rt}P_tf_1>g\right\}}\left(e^{-rt}P_tf_1-e^{-rt}P_tf_0\right) \\ &=& e^{-rt}\chi_{\left\{e^{-rt}P_tf_1>g\right\}}P_t\left(f_1-f_0\right)\end{eqnarray*} which yields the assertion as $P_t$ is an $L^p(\lambda^d)$-contraction (for $p=\infty$ this is immediate and for $p=1$ it follows from the translation-invariance of both $P_t$ and the Lebesgue measure). 
\end{proof}

\begin{rem} Suppose for the moment that $g$ is increasing. Whenever $x\geq -\min_i y_i^{(s)}+z$ componentwise for some $s\in I$ and some $z\in\RR^d$ such that $g(z)\geq 0$ and $g$ is monotonely increasing, then one will have that $P_s(g\vee 0)(x)=P_sg(x)$, thus $\left\{P_s(g\vee 0)=P_sg\right\}$ is not only non-empty, but has even infinite Lebesgue measure (whereas this set is empty if $P_\cdot$ is the semigroup of a diffusion and $g$ can be any continuous function that is not Lebesgue-almost everywhere nonnegative). 

A similar assertion holds for decreasing $g$: Consider any $s\in I$. Whenever $x\leq -\max_i y_i^{(s)}+z$ componentwise for some $z\in\RR^d$ such that $g(z)\leq 0$ and $g$, then one will again have that $P_s(g\vee 0)(x)=P_sg(x)$, thus in this case also $\left\{P_s(g\vee 0)=P_sg\right\}$ has infinite Lebesgue measure. 
\end{rem}

Thus, the statement of the subsequent Lemma does not refer to a null set, let alone the empty set.

\begin{lem} \label{ggeq0onCEt} Assume that either 
\begin{enumerate}
\item $g$ is monotonely decreasing (ie $g(x)\leq g(y)$ whenever $x\geq y$ componentwise for $x,y\in\RR^d$), and
\item $y_k^{(h)}\geq 0$ (componentwise) for some $k\in\{1,\dots,m\}$ (which implies the existence of a $y_i^{(s)}\geq 0$ componentwise for some $i\in\left\{1,\dots,m^\frac{s}{h}\right\}$ for all $s\in I=h\NN_0$), 
\end{enumerate}
or both
\begin{enumerate}
\item $g$ is monotonely increasing, and 
\item $y_k^{(h)}\leq 0$ (componentwise) for some $k\in\{1,\dots,m\}$. 
\end{enumerate}
In both situations, $g$ is nonnegative on $\left\{P_s(g\vee 0)=P_sg\right\}$ for all $s\in I$.
\end{lem}
\begin{proof} On the one hand \begin{eqnarray*}\left\{P_s(g\vee 0)=P_s g\right\}&=& \left\{P_s(g\wedge 0)=0\right\}=\left\{\forall i\in\left\{1,\dots,m^\frac{s}{h}\right\}\quad g\left(\cdot+y_i^{(s)}\right)\geq 0\right\}\\ &=& \bigcap_{i\in\left\{1,\dots,m^\frac{s}{h}\right\}}\left\{  g\left(\cdot+y_i^{(s)}\right)\geq 0\right\}.\end{eqnarray*} 
Now, if $y_k^{(h)}\geq 0$, then $0\leq \frac{s}{h}y_k^{(h)} \in\left\{ y_i^{(s)} \ : \ i\in\left\{1,\dots,m^\frac{s}{h}\right\}\right\}$ for all $s\in I$, and if $y_k^{(s)}\leq 0$, then $0\geq \frac{s}{h}y_k^{(h)} \in\left\{ y_i^{(s)} \ : \ i\in\left\{1,\dots,m^\frac{s}{h}\right\}\right\}$ for all $s\in I$. Thus if $g$ is increasing and $y_k^{(h)}\leq 0$,  or alternatively, $g$ is decreasing and $y_k^{(h)}\geq 0$, then in both cases $g\geq g\left(\cdot + \frac{t}{h}y_k^{(h)}\right)$ and $\frac{t}{h}y_k^{(h)}\in \left\{ y_i^{(t)} \ : \ i\in\left\{1,\dots,m^\frac{t}{h}\right\}\right\}$, therefore  \begin{eqnarray*}\left\{P_s(g\vee 0)=P_s g\right\}&= & \bigcap_{i\in\left\{1,\dots,m^\frac{s}{h}\right\}}\left\{  g\left(\cdot+y_i^{(s)}\right)\geq 0\right\} \\ &=& \left\{ \underbrace{g\left(\cdot+\frac{s}{h}y_k^{(h)}\right)}_{\leq g}\geq 0\right\}\subseteq \left\{ g\geq  0\right\}.\end{eqnarray*} 

\end{proof}

This upper bound on the set $\left\{P_t(g\vee 0) = P_tg\right\}$ (with respect to the partial order of inclusion), can be sharpened under certain conditions on the growth of $\bar f$:

\begin{lem} \label{offEformula} Suppose $g=K-\bar f$ and let there exist a $\gamma_0\geq 1$ (without loss of generality, $\gamma_0\in[1,e^r)$) such that $$P_t\bar f\geq {\gamma_0}^t\bar f$$ for all $t\in(0,T]\cap I$ (which is equivalent to $P_h\bar f\geq {\gamma_0}^h\bar f$). In addition, assume that $g\geq 0$ on the subset $\left\{P_t(g\vee 0)=P_tg\right\}$ of $\RR^d$ for all $t\in(0,T]\cap I$ (by Lemma \ref{ggeq0onCEt}, this is satisfied in particular if $y_k^{(h)}\geq 0$ for some $k\in\{1,\dots,m\}$). Then for all $t\in(0,T]\cap I$, \begin{eqnarray*} \left\{g\geq e^{-rt}P_t(g\vee 0)\right\} &\supseteq&  \left\{P_t(g\vee 0) = P_tg\right\} \end{eqnarray*}
\end{lem}
\begin{proof} Let $t\in (0,T]\cap I$. Due to our assumption of $g\geq 0$ on $\{P_t(g\vee 0)=P_tg\}$, one has \begin{eqnarray} && \left\{g\geq e^{-rt}P_t(g\vee 0)\right\}\cap \left\{P_t(g\vee 0)=P_tg\right\} \nonumber \\ &=& \left\{g\geq e^{-rt}P_tg\right\}\cap \left\{P_t(g\vee 0)=P_tg\right\}  \nonumber \\ &=&\{g\geq 0\}\cap \left\{\left(\id-e^{-rt}P_t\right) g \geq 0\right\}\cap \left\{P_t(g\vee 0)=P_tg\right\} \nonumber  \\ &=& \{g\geq 0\}\cap \left\{\left(1-e^{-rt}\right)K\geq \left(\id-e^{-rt}P_t\right)\bar f\right\} \nonumber \\ &&\cap \left\{P_t(g\vee 0)=P_tg\right\}. \label{lemma1.2eq}\end{eqnarray}
On the other hand \begin{eqnarray*}&&\left\{\left(1-e^{-rt}\right)K\geq \left(\id-e^{-rt}P_t\right)\bar f\right\} = \left\{\left(1-e^{-rt}\right)K\geq \bar f-e^{-rt}\underbrace{P_t\bar f}_{\geq {\gamma_0}^t\bar f}\right\} \\ &\supseteq&   \left\{\left(1-e^{-rt}\right)K\geq \left(1-e^{-rt}{\gamma_0}^t\right)\bar f\right\}, \end{eqnarray*} that is \begin{eqnarray*}&&\left\{\left(1-e^{-rt}\right)K\geq \left(\id-e^{-rt}P_t\right)\bar f\right\} \\ &\supseteq& \left\{\frac{1-e^{-rt}}{1-e^{-rt}{\gamma_0}^t}K \geq  \bar f\right\}, \end{eqnarray*} where we have exploited $\gamma_0<e^r$. Now $\gamma_0\in[1,e^r)$ gives $$\frac{1-e^{-rt}}{1-e^{-rt}{\gamma_0}^t}K\geq K$$ since $K\geq 0$. Combining this estimate with the previous inclusion, one obtains $$\left\{\left(1-e^{-rt}\right)K \geq \left(\id-e^{-rt}P_t\right)\bar f\right\} \supseteq \{K\geq \bar f\}$$ and hence $$\left\{g\geq e^{-rt}P_tg \right\}\cap\{g\geq 0\} =\left\{\left(1-e^{-rt}\right)K\geq \left(\id-e^{-rt}P_t\right)\bar f\right\} \cap \{g\geq 0\} =\{g\geq 0\}.$$ This result, combined with equation (\ref{lemma1.2eq}), yields \begin{eqnarray*} && \left\{g\geq e^{-rt}P_t(g\vee 0)\right\}\cap \left\{P_t(g\vee 0)=P_tg\right\} \\ &=& \{g\geq 0\}\cap \left\{P_t(g\vee 0)=P_tg\right\}. \end{eqnarray*} However, one of our assumptions reads $$\left\{P_t(g\vee 0)=P_tg\right\}\subseteq \{g \geq 0\}$$ whence we conclude \begin{eqnarray*} && \left\{g\geq e^{-rt}P_t(g\vee 0)\right\}\cap \left\{P_t(g\vee 0)=P_tg\right\} \\ &=& \left\{P_t(g\vee 0)=P_tg\right\}. \end{eqnarray*} 
\end{proof}

\begin{rem} Both the assumption of the existence of a constant $\gamma_0\leq e^r$ such that $P_t\bar f\geq {\gamma_0}^t\bar f$ for all $t\in I$, and the assumption of the existence of a constant $\gamma_1>0$ such that $P_t\bar f\leq {\gamma_1}^t\bar f$ for all $t\in I$ are natural: If $(X_t)_{t\in I}$ was a Markov process evolving according to $(P_t)_{t\in I}$, the condition $$\forall t\in I \quad P_t\bar f\geq e^{rt }\bar f$$ simply means that the process $\bar f(X_\cdot)$ is, after discounting, a submartingale. Now, if $X_\cdot$ was a Markovian model for a vector of logarithmic asset prices (a {\em Markov basket} in our terminology) and $\bar f$ would assign to each vector the arithmetic average of the exponentials of its components or the exponential maximum of the components (which equals the maximum of the exponential components) this will certainly hold whenever $P_\cdot$ governs the process $X_\cdot$ under a risk-neutral measure. Furthermore, the assumption $$\forall t\in I \quad P_t\bar f\leq e^{+rt }\bar f$$ (which, being a consequence of the supermartingale property for $e^{-r\cdot}P_\cdot\bar f(x)$ for all $x$, holds whenever $P_\cdot$ is the semigroup of a Markov basket and $\bar f$ assigns to each vector the arithmetic average of the exponential components or the exponential minimum of the components) trivially implies $$\exists \gamma_1>0\quad \forall t\in I\quad P_t\bar f \leq {\gamma_1}^t\bar f$$ and therefore provides us with some economic vindication for assuming the existence of not only $\gamma_0$, but also $\gamma_1$.
\end{rem}

Now we turn to the derivation of an upper bound on the {\em first order difference} $\left(B_{s/2}\right)^{\circ 2}f-B_sf$.

\begin{lem} \label{Bs/2-Bsestimate} Suppose $g:=K-\bar f$ and let there exist a $\gamma_1>0$ such that $$P_t\bar f\leq {\gamma_1}^t\bar f$$ for all $t\in(0,T]\cap I$ (for which in our case of $I=h\NN_0$ with $h>0$ it is sufficient that this estimate holds for $t=h$), and let us assume without loss of generality that this $\gamma_1$ be $\geq e^r$. Then, setting $$R:=K\cdot \sup_{t\in(0,T]\cap I}\frac{{\gamma_1}^t -1}{t},$$ we have found an $R<+\infty$ such that for all $s\in(0,T)\cap (2\cdot I)\subset I$ and measurable $f\geq g\vee 0$, $$\left\|\left(B_{\frac{s}{2}}\right)^{\circ 2}f - B_sf\right\|_{L^\infty(\RR^d)}\leq R\cdot\frac{s}{2}.$$
\end{lem}
\begin{proof} Let $s\in(0,T]\cap (2\cdot I)\subset I$ and consider a measurable $f\geq g\vee 0$. Then by our assumption of $P_t\bar f\leq {\gamma_1}^t\bar f$ for $\gamma_1>0$, we firstly have (inserting $\frac{s}{2}$ for $t$)
\begin{eqnarray*}&& g-e^{-r\frac{s}{2}}P_{\frac{s}{2}}g=K-\bar f-e^{-r\frac{s}{2}}\left(K-P_{\frac{s}{2}}\bar f\right)\\ &=& K\left(1-e^{-r\frac{s}{2}}\right)+e^{-r\frac{s}{2}}P_{\frac{s}{2}}\bar f-\bar f\\ &\leq& K\left(1-e^{-r\frac{s}{2}}\right)+ \left(e^{-r\frac{s}{2}}{\gamma_1}^\frac{s}{2}-1\right)\cdot\bar f \end{eqnarray*}
and therefore (using $f\geq g \vee 0$ and $\gamma_1\geq e^r$ as well as the monotonicity of $P_{\frac{s}{2}}$), 
\begin{eqnarray} 0&\leq& \chi_{\left\{g\geq e^{-r\frac{s}{2}}P_{\frac{s}{2}}f\right\}}\cdot\left(g-e^{-r\frac{s}{2}}P_{\frac{s}{2}}f\right) \nonumber\\ &\leq& \chi_{\left\{g\geq e^{-r\frac{s}{2}}P_{\frac{s}{2}}f\right\}}\cdot\left(g-e^{-r\frac{s}{2}}P_{\frac{s}{2}}g\right) \nonumber\\ &\leq& \chi_{\left\{g\geq 0\right\}}\cdot\left(K\left(1-e^{-r\frac{s}{2}}\right)+ \left(e^{-r\frac{s}{2}}{\gamma_1}^\frac{s}{2}-1\right)\cdot\bar f \right) \nonumber\\ &\leq& \chi_{\left\{\bar f\leq K\right\}}\cdot\left(K\left(1-e^{-r\frac{s}{2}}\right)+ \left(e^{-r\frac{s}{2}}{\gamma_1}^\frac{s}{2}-1\right)\cdot K \right)\nonumber \\ &\leq& K\cdot\left(\left(1-e^{-r\frac{s}{2}}\right)+ \left(e^{-r\frac{s}{2}}{\gamma_1}^\frac{s}{2}-1\right) \right)\nonumber \\ & = & K\cdot e^{-r\frac{s}{2}}\cdot \left({\gamma_1}^\frac{s}{2}-1\right) \label{Rbound} \end{eqnarray} 
Now, $$\sup_{t\in(0,T]\cap I}\frac{{\gamma_1}^t -1}{t}<+\infty$$ since $t\mapsto {\gamma_1}^t-1$ is right-differentiable in zero with derivative $\ln \gamma_1$. Therefore $$R=K\cdot \sup_{t\in(0,T]\cap I}\frac{{\gamma_1}^t -1}{t}<+\infty.$$ Via estimate (\ref{Rbound}), we arrive at \begin{equation}0\leq \chi_{\left\{g\geq e^{-r\frac{s}{2}}P_{\frac{s}{2}}f\right\}}\left(g-e^{-r\frac{s}{2}}P_{\frac{s}{2}}f\right)\leq R\cdot \frac{s}{2}.\label{Chapter8,1.3a}\end{equation} But $$\chi_{\left\{g\geq e^{-r\frac{s}{2}}P_{\frac{s}{2}}f\right\}}\left(g-e^{-r\frac{s}{2}}P_{\frac{s}{2}}f\right)=\left(e^{-r\frac{s}{2}}P_{\frac{s}{2}}f\right)\vee g - e^{-r\frac{s}{2}}P_{\frac{s}{2}}f,$$ and -- in combination with the linearity of $P_{\frac{s}{2}}$ and the Chapman-Kolmogorov equation -- this implies \begin{eqnarray*}&& e^{-r\frac{s}{2}}P_{\frac{s}{2}}\left(\chi_{\left\{g\geq e^{-r\frac{s}{2}}P_{\frac{s}{2}}f\right\}}\left(g-e^{-r\frac{s}{2}}P_{\frac{s}{2}}f\right)\right)\\ &=&e^{-r\frac{s}{2}}P_{\frac{s}{2}}\left(\left(e^{-r\frac{s}{2}}P_{\frac{s}{2}}f\right)\vee g\right) - e^{-rs}P_{{s}}f\end{eqnarray*} hence by equation (\ref{Chapter8,1.3a})\begin{eqnarray}&& R\cdot\frac{s}{2}\nonumber\\ &\geq& e^{-r\frac{s}{2}}P_{\frac{s}{2}}\left(\chi_{\left\{g\geq e^{-r\frac{s}{2}}P_{\frac{s}{2}}f\right\}}\left(g-e^{-r\frac{s}{2}}P_{\frac{s}{2}}f\right)\right)\\ &\geq&e^{-r\frac{s}{2}}P_{\frac{s}{2}}\left(\left(e^{-r\frac{s}{2}}P_{\frac{s}{2}}f\right)\vee g\right) - e^{-rs}P_{{s}}f\vee g \label{Rs/2eqn}\end{eqnarray} (where in (\ref{Rs/2eqn}) we have exploited the fact that $P_\frac{s}{2}$ is an $L^\infty$-contraction).

Now, again by the Chapman-Kolmogorov equation and the monotonicity of $P_t$ for any $t\in I$, we have $$\left(B_{\frac{s}{2}}\right)^{\circ 2}f\geq B_s f$$ and therefore (due to the estimate $B_tf_0\geq g\vee 0$ which holds for arbitrary $t\in I$ and $f_0\geq 0$) \begin{eqnarray*} \left\{\left(B_{\frac{s}{2}}\right)^{\circ 2}f=g\right\}&=& \left\{B_sf\leq \left(B_{\frac{s}{2}}\right)^{\circ 2}f=g\right\} \\ &=&  \left\{g\vee 0 \leq B_sf\leq \left(B_{\frac{s}{2}}\right)^{\circ 2}f=g\right\} \\ &=& \left\{B_sf=g\right\}\cap  \left\{\left(B_{\frac{s}{2}}\right)^{\circ 2}f=g\right\}\\ &\subseteq&  \left\{\left(B_{\frac{s}{2}}\right)^{\circ 2}f-B_sf=0\right\} \end{eqnarray*} 
But $$\left(B_{\frac{s}{2}}\right)^{\circ 2}f\geq g\vee 0\geq g,$$ thus the last inclusion yields \begin{eqnarray*}0&\leq& \left(B_{\frac{s}{2}}\right)^{\circ 2}f-B_sf\\ &=&\left(\left(B_{\frac{s}{2}}\right)^{\circ 2}f-B_sf\right)\chi_{\left\{ \left(B_{\frac{s}{2}}\right)^{\circ 2}f\geq g\right\}}\\ &=& e^{-r\frac{s}{2}}\left(P_{\frac{s}{2}}\left(e^{-r\frac{s}{2}}P_{\frac{s}{2}}f\vee g\right) - \left(e^{-rs}P_{{s}}f\vee g\right)\right)\chi_{\left\{ \left(B_{\frac{s}{2}}\right)^{\circ 2}f\geq g\right\}}\\ &\leq& R\cdot\frac{s}{2}\end{eqnarray*} where the last line has used the estimate (\ref{Rs/2eqn}) derived previously.
\end{proof}

A bit further on, in Lemma \ref{lowerboundoffE}, we will establish lower bounds for $\left\|\left(B_{\frac{s}{2}}\right)^{\circ 2}f - B_sf\right\|_{L^\infty(\RR^d)}$ that are uniform in $f\geq g\vee 0$ and only linear in $s$. Therefore it is impossible to get from estimates on $\left\|\left(B_{\frac{s}{2}}\right)^{\circ 2}f - B_sf\right\|_{L^\infty(\RR^d)}$ to estimates on $\left\|\left(B_{s2^{-(k+1)}}\right)^{\circ 2^{k+1}}f - \left(B_{s{2^-k}}\right)^{\circ 2^k}f\right\|_{L^\infty(\RR^d)}$ in the vein of Lemma \ref{from_k=1_to_kinN} (where one uses estimates on $\left(B_{\frac{s}{2}}\right)^{\circ 2}f - B_sf$ that are of higher than linear order in $s$).

We can draw from the proof of Lemma \ref{Bs/2-Bsestimate} the following Corollary which can be used for instance to prove estimates on the difference $\left(B_{\frac{s}{2}}\right)^{\circ 2}f-B_sf$ with respect to the $L^\infty \left\{P_t(g\vee 0)>P_tg\right\}$- and $L^1\left\{P_t(g\vee 0)>P_tg\right\}$-norms.

\begin{cor}\label{Bs/2-Bsequation} For arbitrary $g:\RR^d\rightarrow\RR$ and all $f\geq g\vee 0$, \begin{eqnarray*}0 &\leq& \left(B_{\frac{s}{2}}\right)^{\circ 2}f-B_sf\\ &=& e^{-r\frac{s}{2}}P_{\frac{s}{2}}\left(\left(e^{-r\frac{s}{2}}P_{\frac{s}{2}}f\right)\vee g - e^{-r\frac{s}{2}}P_{\frac{s}{2}}f\right)\chi_{\left\{ \left(B_{\frac{s}{2}}\right)^{\circ 2}f\geq g\right\}}\\&\leq& e^{-r\frac{s}{2}}P_{\frac{s}{2}}\left(\chi_{\left\{g\geq e^{-r\frac{s}{2}}P_{\frac{s}{2}}(g\vee 0)\right\}}\cdot\left(g-e^{-r\frac{s}{2}}P_{\frac{s}{2}}(g\vee 0)\right) \right)\\ &=& \left(B_{\frac{s}{2}}\right)^{\circ 2}(g\vee 0)-B_s(g\vee 0).\end{eqnarray*}
\end{cor}

\begin{Def}$$\forall t\in I\quad E^t = \left\{P_t(g\vee 0)>P_tg\right\}.$$\end{Def}

\begin{rem}[Properties of $E^\cdot$] \label{Etproperties} Let $g=K-\bar f$. Then equivalent expressions for $E^t$, $t\in I$, are:
\begin{eqnarray*} E^t &=& \left\{P_t(g\vee 0)>P_tg\right\} = \complement \left\{P_t(g\vee 0)=P_tg\right\} \\ &=& \left\{P_t(g\wedge 0)<0\right\}  \\ &=& \left\{\exists i\in\{1,\dots,m^{\frac{t}{h}}\}\quad g\left(\cdot+y_i^{(t)}\right)<0\right\}\\ &=& \bigcup_{i\in\{1,\dots,m^{\frac{t}{h}}\}}\left\{K<\bar f\left(\cdot+y_i^{(t)}\right)\right\} .\end{eqnarray*} These formulae for $E^\cdot$ imply, since $\bar f$ is monotonely increasing, that $E^t$ is north-east connected, ie $E^t\subseteq E^t+a$ for all $a\leq 0$, for all $t\in I$. 

Assume now that one has for all $j\in\left\{1,\dots, m\right\}$ an index $k=k(j)\in\left\{1,\dots, m\right\}$ such that $y_j^{(h)}-y_{k(j)}^{(h)}\geq 0$ componentwise (this assumption is satisfied for instance if the set $\left\{y_1^{(h)}, \dots,y_m^{(h)} \right\}$ can be written as the sum of a reflection symmetric subset of $\RR^d$ and a componentwise nonpositive vector). Then for any $t\in I$ and $i\in\left\{1,\dots, m^\frac{t}{h}\right\}$ there exists an index $k=k(i)\in\left\{1,\dots, m^\frac{t}{h}\right\}$ such that $y_i^{(t)}-y_{k(i)}^{(t)}\geq 0$ componentwise. Therefore the north-east connectedness of the $E^t$'s entails for all $t\in I$ and $i$, \begin{eqnarray*}E^t+y_i^{(t)}&\subseteq& E^t\underbrace{-y_i^{(t)}+y_{k(i)}^{(t)}}_{\leq 0}+y_i^{(t)}=E^t+y_{k(i)}^{(t)}\\ &=&\bigcup_{j\in\left\{1,\dots,m^\frac{t}{h}\right\}}\left\{K<\bar f\left(\cdot+y_j^{(t)}+y_{k(i)}^{(t)}\right)\right\} \\ &\subseteq&\bigcup_{j_0,j_1\in\left\{1,\dots,m^\frac{t}{h}\right\}}\left\{K<\bar f\left(\cdot+y_{j_0}^{(t)}+y_{j_1}^{(t)}\right)\right\}\\ &=&\bigcup_{\ell\in\left\{1,\dots,m^\frac{2t}{h}\right\}}\left\{K<\bar f\left(\cdot+y_\ell^{(2t)}\right)\right\} \\ &=& E^{2t}\end{eqnarray*} where for the penultimate line we have used equation \eqref{chapman}, of course. Therefore $$\chi_{E^t}\left(\cdot-y_i^{(t)}\right)=\chi_{E^t+y_i^{(t)}}\leq\chi_{E^{2t}}.$$ Also, if there exists an $i_0\in\{1,\dots,m\}$ such that $$\forall j\in\{1,\dots, d\}\quad \left(y_{i_0}^{(h)}\right)_j\geq 0 $$ one has -- due to the monotonicity of $\bar f$ in each coordinate -- first of all $\bar f\left(\cdot +y_{i_0}^{(h)}\right)\geq \bar f$ and thence for all $n\in\NN$ the inclusion \begin{eqnarray*} E^{nh} &=& \bigcup_{i\in\{1,\dots,m^n\}}\left\{K<\bar f\left(\cdot+y_i^{(nh)}\right)\right\}  \\ &=& \bigcup_{i_1,\dots,i_n\in\{1,\dots,m\}}\left\{K<\bar f\left(\cdot+y_{i_1}^{(h)}+\dots +y_{i_n}^{(h)}\right)\right\}\\ &\supseteq& \bigcup_{i_1,\dots,i_{n-1}\in\{1,\dots,m\}}\left\{K<\bar f\left(\cdot+y_{i_0}^{(h)}+y_{i_1}^{(h)}+\cdots +y_{i_{n-1}}^{(h)}\right)\right\} \\ &\supseteq& \bigcup_{i_1,\dots,i_{n-1}\in\{1,\dots,m\}}\left\{K<\bar f\left(\cdot+y_{i_1}^{(h)}+\cdots +y_{i_{n-1}}^{(h)}\right)\right\}  \\ &=& E^{(n-1)h} \end{eqnarray*} This means $$E^s\uparrow \text{ as }s\uparrow\infty\text{ in }I$$ and for all $T\in[h,+\infty]$, $$\bigcap_{s\in (0,T]\cap I}E^s = E^h.$$
\end{rem}

The reason for $E^h$ not being the whole space is, as was pointed our earlier, that the measure $B\mapsto P_t\chi_B$ on the Borel $\sigma$-algebra of $\RR$ has compact support. 

Interpreting $g$ as defining a logarithmic payoff function (eg $g=K-\exp$, $d=1$ in case of a vanilla one-dimensional put) and $P$ as a Markov chain that models the stochastic evolution of the logarithmic prices of assets in a given portfolio, the set $E^t$, for $t\in I$, consists of all those vectors of logarithmic start prices where the probability of exercising the option at time $t$ is strictly positive.

\begin{rem} Assume $\bar f$ is not strictly less than $K$, say $\bar f(z)\geq K$ for some $z\in\RR^d$. We can use the property of $\bar f$ being monotonely increasing in each component to see, via Remark \ref{Etproperties} that 
\begin{eqnarray*}\complement E^h&=&\left\{P_h(g\vee 0)=P_hg\right\} \\ &=&\complement\left\{\exists i\in\{1,\dots,m\}\quad g\left(\cdot+y_i^{(h)}\right)<0\right\} \\ &=&\bigcap_{i=1}^m\left\{g\left(\cdot+y_i^{(h)}\right)\geq 0\right\} \\ &=&\bigcap_{i=1}^m\left\{\bar f\left(\cdot+y_i^{(h)}\right) \leq K \leq \bar f(z)\right\} \\ &\supseteq &\bigcap_{i=1}^m\left\{\forall j\in\{1,\dots,d\}\quad \left(\cdot+y_i^{(h)}\right)_j \leq z_j \right\} \\ &= &\bigcap_{i=1}^m\left\{\forall j\in\{1,\dots,d\}\quad \left(\cdot\right)_j \leq z_j -\left(y_i^{(h)}\right)_j \right\} \\ &= & \left\{\forall j\in\{1,\dots,d\}\quad \left(\cdot\right)_j \leq z_j -\max_{i\in\{1,\dots,m\}}\left(y_i^{(h)}\right)_j \right\} \\ &= & \bigotimes_{j=1}^d\left(-\infty,z_j- \max_{i\in\{1,\dots,m\}}\left(y_i^{(h)}\right)_j \right] ,\end{eqnarray*} where the set in the last line has infinite Lebesgue measure. 

Thus, $\lambda^d\left[\complement E^h\right]=+ \infty$ whenever $\bar f< K$ does not hold everywhere.

\end{rem}

Next, we can combine the preceding Corollary \ref{Bs/2-Bsequation} (which followed from the proof of Lemma \ref{Bs/2-Bsestimate}) with Lemmas \ref{ggeq0onCEt} and \ref{offEformula} to prove the previously mentioned lower bound on the $L^\infty$-norm of $\left(B_{\frac{t}{2}}\right)^{\circ t}f - B_tf$ on the complement of the set $E^\frac{t}{2}$ introduced previously.

\begin{lem} \label{lowerboundoffE} Let $g=K-\bar f$. Suppose there is a $\gamma_0>1$ (without loss of generality, $\gamma_0\in\left(1,e^r\right]$) and a $\gamma_1>0$ such that $${\gamma_1}^t\bar f\geq P_t\bar f\geq {\gamma_0}^t\bar f$$ for all $t\in(0,T]\cap I$ (where  $I=h\NN_0$ with $h>0$ whence it is sufficient that this estimate holds for $t=h$). 
Assume furthermore that $y_{i_0}^{(h)}\leq 0$ for some $i_0\in\{1,\dots,m\}$ and $y_{i_1}^{(h)}\geq 0$ for some $i_1\in\{1,\dots,m\}$, implying that $g\geq 0$ on the set $\left\{P_t(g\vee 0)=P_tg\right\}\neq\emptyset$. Then for all $\varepsilon_1>0$ there is an $\varepsilon_0<T$ independent of $h<T$ such that for all $t\in 2\cdot\left( \left(0,\varepsilon_0\right)\cap I \right)$ and $A\supset \complement E^\frac{t}{2}$ (with positive Lebesgue measure), \begin{eqnarray*}&&\sup_{f\geq g\vee 0}\left\|\left(B_{\frac{t}{2}}\right)^{\circ 2}f - B_tf\right\|_{L^\infty\left(A \right)}\geq \left\|\left(B_{\frac{t}{2}}\right)^{\circ 2}(g\vee 0) - B_t(g\vee 0)\right\|_{L^\infty\left(A \right)} \\ &\geq& \left(\min_{i\in\left\{1,\dots,m\right\}} \alpha^{\left(h\right)}\right)^\frac{T}{h} e^{-r\frac{t}{2}}K\left(\ln\gamma_0-\varepsilon_1 \right)\cdot \frac{t}{2},\end{eqnarray*} as well as \begin{eqnarray*}&&\sup_{f\geq g\vee 0}\left\|\left(B_{\frac{t}{2}}\right)^{\circ 2}f - B_tf\right\|_{L^1\left(A \right)}\geq \left\|\left(B_{\frac{t}{2}}\right)^{\circ 2}(g\vee 0) - B_t(g\vee 0)\right\|_{L^1\left(A \right)} \\ &\geq& \left(\min_{i\in\left\{1,\dots,m\right\}} \alpha^{\left(h\right)}\right)^\frac{T}{h} \lambda^d\left[\left\{P_{\frac{t}{2}}(g\vee 0)=P_{\frac{t}{2}} g\right\}\right]\cdot e^{-r\frac{t}{2}}K\left(\ln\gamma_0-\varepsilon_1 \right)\cdot \frac{t}{2} \end{eqnarray*} (the left hand side, following the usual convention, being $+\infty$ if $\lambda^d\left[\left\{P_{\frac{t}{2}}(g\vee 0)=P_{\frac{t}{2}} g\right\}\right]=+\infty$, $m \geq 1$ and $\varepsilon_1<\ln\gamma_0$).
\end{lem}

\begin{proof} Let us first remark that, due to Corollary \ref{Bs/2-Bsequation}, we have \begin{eqnarray} \nonumber &&\sup_{f\geq g\vee 0} \left\|\left(B_{\frac{t}{2}}\right)^2f-B_{t}f\right\|_{L^\infty\left(A\right)} \\ \nonumber & = & \left\|\left(B_{\frac{t}{2}}\right)^2(g\vee 0)-B_{t}(g\vee 0)\right\|_{L^\infty\left(A\right)}\\ &=&  \left\|\chi_A e^{-r\frac{t}{2}}P_{\frac{t}{2}}\left(\chi_{\left\{g\geq e^{-r\frac{t}{2}}P_{\frac{t}{2}}(g\vee 0)\right\}}\cdot\left(g-e^{-r\frac{t}{2}}P_{\frac{t}{2}}(g\vee 0)\right)\right) \right\|_{L^\infty\left(A\right)} \label{lowerboundoffEeq1}\end{eqnarray} 
as well as 
\begin{eqnarray} &&\sup_{f\geq g\vee 0} \left\|\left(B_{\frac{t}{2}}\right)^2f-B_{t}f\right\|_{L^1\left(A\right)} \nonumber\\ & = & \left\|\left(B_{\frac{t}{2}}\right)^2(g\vee 0)-B_{t}(g\vee 0)\right\|_{L^\infty\left(A\right)}\nonumber\\ &=&  \left\|e^{-r\frac{t}{2}}P_{\frac{t}{2}}\left(\chi_{\left\{g\geq e^{-r\frac{t}{2}}P_{\frac{t}{2}}(g\vee 0)\right\}}\cdot\left(g-e^{-r\frac{t}{2}}P_{\frac{t}{2}}(g\vee 0)\right)\right) \right\|_{L^1\left(A\right)}\nonumber\\ &=&  \left\|\left(\chi_A\cdot e^{-r\frac{t}{2}}P_{\frac{t}{2}}\left(\chi_{\left\{g\geq e^{-r\frac{t}{2}}P_{\frac{t}{2}}(g\vee 0)\right\}}\cdot\left(g-e^{-r\frac{t}{2}}P_{\frac{t}{2}}(g\vee 0)\right)\right) \right)\left(\cdot+\frac{t}{2h}y_{i_0}^{(h)}\right)\right\|_{L^1\left(\RR^d\right)}\label{chapter8a,4*}\end{eqnarray} 
for all $i_0\in \left\{1,\dots,m^\frac{t}{2h}\right\}$ (using the translation invariance of $\int_{\RR^d}\cdot\lambda^d$).

Next let us note that by our assumption of $y_{i_0}^{(h)}\leq 0$ componentwise, combined with the north-east connectedness of $E^s$, we have $\frac{s}{h} y_{i_0}^{(h)}+E^s\subseteq E^s$ for any $s$, hence $$\forall s\in I\quad \frac{s}{h} y_{i_0}^{(h)}+\complement E^s\supseteq \complement E^s.$$ Therefore we may conclude that for all $t\in 2\cdot I$ and $i_0\in \left\{1,\dots,m^\frac{t}{2h}\right\}$, 
\begin{eqnarray}&&\chi_{A} \cdot e^{-r\frac{t}{2}}P_{\frac{t}{2}}\left(\chi_{\left\{g\geq e^{-r\frac{t}{2}}P_{\frac{t}{2}}(g\vee 0)\right\}}\cdot\left(g-e^{-r\frac{t}{2}}P_{\frac{t}{2}}(g\vee 0)\right)\right) \nonumber \\ &\geq&  \chi_{\complement E^\frac{t}{2}} \cdot e^{-r\frac{t}{2}}P_{\frac{t}{2}}\left(\chi_{\left\{g\geq e^{-r\frac{t}{2}}P_{\frac{t}{2}}(g\vee 0)\right\}}\cdot\left(g-e^{-r\frac{t}{2}}P_{\frac{t}{2}}(g\vee 0)\right)\right) \nonumber \\ &\geq& \chi_{\complement E^\frac{t}{2}}\cdot \left(\chi_{\left\{g\geq e^{-r\frac{t}{2}}P_{\frac{t}{2}}(g\vee 0)\right\}}\cdot\left(g-e^{-r\frac{t}{2}}P_{\frac{t}{2}}(g\vee 0)\right)\right)\left(\cdot+\frac{t}{2h}y_{i_0}^{(h)}\right)\nonumber \\&& \cdot  \min_{i\in\left\{1,\dots,m^\frac{t}{2h}\right\}}{\alpha^{\left(\frac{t}{2}\right)}}  \nonumber \\ &\geq& \left(\chi_{\left\{g\geq e^{-r\frac{t}{2}}P_{\frac{t}{2}}(g\vee 0)\right\}}\cdot\left(g-e^{-r\frac{t}{2}}P_{\frac{t}{2}}(g\vee 0)\right)\chi_{\complement E^\frac{t}{2}+\frac{t}{2h} y_{i_0}^{(h)}} \right)\left(\cdot+\frac{t}{2h}y_{i_0}^{(h)}\right)\nonumber \\&& \cdot  \left(\min_{i\in\left\{1,\dots,m\right\}} \alpha^{\left(h\right)}\right)^\frac{t}{2h}\nonumber  \\ &\geq& \left(\chi_{\left\{g\geq e^{-r\frac{t}{2}}P_{\frac{t}{2}}(g\vee 0)\right\}}\cdot\left(g-e^{-r\frac{t}{2}}P_{\frac{t}{2}}(g\vee 0)\right)\chi_{\complement E^\frac{t}{2}}\right)\left(\cdot+\frac{t}{2h}y_{i_0}^{(h)}\right)\nonumber \\&& \cdot  \left(\min_{i\in\left\{1,\dots,m\right\}} \alpha^{\left(h\right)}\right)^\frac{t}{2h}\label{chapter8a,1.6} \\ (& \geq& 0)\nonumber \end{eqnarray}

Now, off $E^s$ one has due to Lemma \ref{offEformula} (which may be applied thanks to our assumption of $y_{i_1}^{(h)}\geq 0$) and Lemma \ref{ggeq0onCEt} the following situation: \begin{eqnarray}&&\chi_{\left\{g\geq e^{-rs}P_{s}(g\vee 0)\right\}}\cdot\left(g-e^{-r{s}}P_{{s}}(g\vee 0)\right) \nonumber \\ &=& \chi_{\left\{P_s(g\vee 0)=P_sg\right\}\cap\{g\geq 0\}}\cdot\left(g-e^{-r{s}}P_{{s}}(g\vee 0)\right) \quad\text{ on }\complement E^s\nonumber  \\ &=& \chi_{\left\{P_s(g\vee 0)=P_sg\right\}\cap\{g\geq 0\}}\cdot\left(g-e^{-r{s}}P_{{s}}g\right) \quad\text{ on }\complement E^s \nonumber  \\ &=& \chi_{\left\{P_s(g\vee 0)=P_sg\right\}\cap\{g\geq 0\}}\cdot\left(K-\bar f-e^{-r{s}}K+e^{-rs}P_s\bar f\right) \quad\text{ on }\complement E^s \nonumber \\ &\geq & \chi_{\left\{P_s(g\vee 0)=P_sg\right\}\cap\{g\geq 0\}}\cdot\left(K-\bar f-e^{-r{s}}K+e^{-rs}{\gamma_0}^s\bar f\right) \quad\text{ on }\complement E^s \label{lowerboundoffEest} \end{eqnarray}
However, one can also perform the calculation 
\begin{eqnarray} && \chi_{\{g\geq 0\}}\left(K-\bar f-e^{-r{s}}K+{\gamma_0}^s e^{-rs}\bar f\right) \nonumber \\ &=& \chi_{\{K-\bar f\geq 0\}}\left(\left(K-\bar f\right)\left(1-{\gamma_0}^se^{-r{s}}\right) +Ke^{-rs}\left({\gamma_0}^s-1\right) \right)\label{lowerboundoffEpenult} \\ &\geq & Ke^{-rs}\left({\gamma_0}^s-1\right) \label{lowerboundoffEfin}\end{eqnarray} (where we have used the assumption $\gamma_0\leq e^r$ to get from (\ref{lowerboundoffEpenult}) to (\ref{lowerboundoffEfin})). Combining estimates (\ref{lowerboundoffEfin}) and (\ref{lowerboundoffEest}), we arrive at \begin{eqnarray*} &&\chi_{\left\{g\geq e^{-rs}P_{s}(g\vee 0)\right\}}\left(g-e^{-r{s}}P_{{s}}(g\vee 0)\right) \\ &\geq& Ke^{-rs}\left({\gamma_0}^s-1\right)  \quad\text{ on }\complement E^s \\ &\geq& Ke^{-rs}\left(\ln{\gamma_0}-\varepsilon_1\right)\cdot s \quad\text{ on }\complement E^s \end{eqnarray*} for every $s<\varepsilon_0$ for some $\varepsilon_0>0$ dependent on $\varepsilon_1>0$ and finally (using estimate (\ref{chapter8a,1.6}), $\min_{i\in\left\{1,\dots,m\right\}} \alpha^{\left(h\right)}\leq 1$ and $T\geq t$) \begin{eqnarray*} && \chi_{A} \cdot e^{-r\frac{t}{2}}P_{\frac{t}{2}}\left(\chi_{\left\{g\geq e^{-r\frac{t}{2}}P_{\frac{t}{2}}(g\vee 0)\right\}}\cdot\left(g-e^{-r\frac{t}{2}}P_{\frac{t}{2}}(g\vee 0)\right)\right) \\ &\geq &   \left(\min_{i\in\left\{1,\dots,m\right\}} \alpha^{\left(h\right)}\right)^\frac{t}{2h}\cdot \chi_{\complement E^\frac{t}{2}}\left(\cdot+\frac{t}{2h}y_{i_0}^{(h)}\right) \cdot Ke^{-r{\frac{t}{2}}}\left(\ln{\gamma_0}-\varepsilon_1\right)\cdot {\frac{t}{2}} \\ & \geq & \left(\min_{i\in\left\{1,\dots,m\right\}} \alpha^{\left(h\right)}\right)^\frac{T}{2h}\cdot \chi_{\complement E^\frac{t}{2} -\frac{t}{2h}y_{i_0}^{(h)}}\cdot Ke^{-r{\frac{t}{2}}}\left(\ln{\gamma_0}-\varepsilon_1\right)\cdot {\frac{t}{2}}\end{eqnarray*} for all $t\in 2\cdot \left(I\cap(0,\varepsilon_0)\right)$ and $i_0\in \left\{1,\dots,m^\frac{t}{2h}\right\}$.

This yields -- due to the translation-invariance of the Lebesgue measure (which gave us estimate (\ref{chapter8a,4*})) -- the first line of the Lemma's $L^1$ norm estimate. It also implies the $L^\infty$ norm estimate of the Lemma since for all $s\in I$ (in particular for $s=\frac{t}{2}$), \begin{eqnarray*}\lambda^d\left[\complement E^{s} -\frac{s}{h}y_{i_0}^{(h)}\right]&=&\lambda^d\left[\complement E^\frac{s}{h}\right] \\ &= &\lambda^d\left\{P_s(g\vee 0)=P_s g\right\}\\ &\geq &\lambda^d\left\{ g\left(\cdot + \frac{s}{h}\max_{i\in\{1,\dots, m\}} y_i^{(h)}\right)\geq 0\right\}>0\end{eqnarray*} (a consequence of the monotonicity of $g$), and therefore \begin{eqnarray*}&&\left\|\left(\min_{i\in\left\{1,\dots,m\right\}} \alpha^{\left(h\right)}\right)^\frac{T}{2h}\cdot \chi_{\complement E^\frac{t}{2}-\frac{t}{2h}y_{i_0}^{(h)}}\cdot Ke^{-r{\frac{t}{2}}}\left(\ln{\gamma_0}-\varepsilon_1\right)\cdot {\frac{t}{2}}\right\|_{L^\infty(\RR^d)}\\&=&\left(\min_{i\in\left\{1,\dots,m\right\}} \alpha^{\left(h\right)}\right)^\frac{T}{2h}\cdot  Ke^{-r{\frac{t}{2}}}\left(\ln{\gamma_0}-\varepsilon_1\right)\cdot {\frac{t}{2}}.\end{eqnarray*}
\end{proof}


\section{Finer bounds on the first order differences on $E^\cdot$}

In Lemma \ref{lowerboundoffE} we have established a lower bound on the $L^1$ and $L^\infty$-norms of $\left(B_{\frac{t}{2}}\right)^{\circ t}f - B_tf$ on the set $\complement \cap_{t\in I}E^t$ that are linear in $t$. We will now try and establish $L^1$-estimates for $\left(B_{\frac{t}{2}}\right)^{\circ t}f - B_tf$ on $\cap_{t\in I}E^t$ (which under weak assumptions equals $E^h$) that are of higher order in $t$.

Keeping Corollary \ref{Bs/2-Bsequation} in mind and aiming at an $L^1\left(\cap_{t\in I}E^t\right)$ of Lemma \ref{Bs/2-Bsestimate}, the first step will consist in proving

\begin{lem} \label{estimateonE} Let again $g= K-\bar f $ and suppose there is a $\gamma_1>0$ such that $$P_t\bar f\leq {\gamma_1}^t\bar f$$ for all $t\in(0,T]\cap I$. Let us define $$\tilde D:=\chi_{(0,e^r)}\left(\gamma_1\right)\inf_{\bigcup_{t\in(0,T]\cap I}E^t}\bar f +\chi_{(e^r,+\infty)}\left(\gamma_1\right)\sup_{\bigcup_{t\in(0,T]\cap I}E^t}\bar f\geq 0.$$
Then there is a constant $C_0\in \RR$ given by $$C_0:=K\left(\sup_{s\in(0,T]\cap I}\frac{1-e^{-rs}}{s}- r\right)+\tilde D\cdot \left(\sup_{s\in(0,T]\cap I} \frac{{\gamma_1}^se^{-rs}-1}{s}-\ln\gamma_1+ r\right)$$ such that for all $s\in(0,T]\cap I$ and measurable $A$, \begin{eqnarray*}&& \left\|\chi_{\left\{g> e^{-rs}P_{s}(g\vee 0)\right\}}\left(g-e^{-r{s}}P_{{s}}(g\vee 0)\right) \right\|_{L^1\left(\bigcap_{s\in (0,T]\cap I}E^s\cap A\right)} \\ &\leq &\left\|\chi_{\left\{g> e^{-rs}P_{s}(g\vee 0)\right\}}\left(g-e^{-r{s}}P_{{s}}(g\vee 0)\right) \right\|_{L^1(E^s\cap A)} \\ &\leq& \lambda^d\left[\left\{e^{rs}g>P_s(g\vee 0)>P_s g\right\}\cap A\right] \cdot \left(\left(\ln\gamma_1-r\right)\tilde D+rK+ C_0 \right)\cdot s .\end{eqnarray*}
If moreover, there exists an $i_0\in\{1,\dots,m\}$ such that for any $j\in\{1,\dots, d\}$, $\left(y_{i_0}^{(h)}\right)_j\geq 0 $ for all $j\in\{1,\dots, d\}$, then we even have for all measurable $A$ \begin{eqnarray*}&& \left\|\chi_{\left\{g> e^{-rs}P_{s}(g\vee 0)\right\}}\left(g-e^{-r{s}}P_{{s}}(g\vee 0)\right) \right\|_{L^1\left(E^h\cap A\right)}  \\ &\leq& \lambda^d\left[\left\{e^{rs}g>P_s(g\vee 0)>P_s g\right\}\cap A\right] \cdot \left(\left(\ln\gamma_1-r\right)\tilde D+rK+ C_0 \right)\cdot s .\end{eqnarray*}
\end{lem}
\begin{proof} For all $s\in(0,T]\cap I$, the following estimates hold on $E^s$: \begin{eqnarray*}0&\leq& \chi_{\left\{g> e^{-rs}P_{s}(g\vee 0)\right\}}\left(g-e^{-r{s}}P_{{s}}(g\vee 0)\right) \\ &=& \chi_{\left\{e^{rs}g>P_s(g\vee 0)>P_sg\right\}}\left(g-e^{-r{s}}P_{{s}}(g\vee 0)\right) \quad\text{ on }E^s  \\ &\leq& \chi_{\left\{e^{rs}g>P_s(g\vee 0)>P_sg\right\}}\left(g-e^{-r{s}}P_{{s}}g\right) \\ &\leq& \chi_{\left\{e^{rs}g>P_s(g\vee 0)>P_sg\right\}}\left(g-e^{-r{s}}P_{{s}}g\right) \\ &=& \chi_{\left\{e^{rs}g>P_s(g\vee 0)>P_sg\right\}}\left(K-\bar f-e^{-rs}K+e^{-rs} P_s\bar f\right) \\ &\leq& \chi_{\left\{e^{rs}g>P_s(g\vee 0)>P_sg\right\}}\left(K\left(1-e^{-rs}\right)+\left({\gamma_1}^se^{-rs}-1\right)\bar f\right) \\ &\leq& \chi_{\left\{e^{rs}g>P_s(g\vee 0)>P_sg\right\}}\left(K\left(1-e^{-rs}\right)+\left({\gamma_1}^se^{-rs}-1\right) \tilde D\right) \\ &\leq& \chi_{\left\{e^{rs}g>P_s(g\vee 0)>P_sg\right\}}\left(rKs+\left(\ln{\gamma_1}-r\right)\tilde D\cdot s +C\cdot s \right)\end{eqnarray*} for some real constant $C>0$ that can be bounded by $$C\leq K\cdot\left(\sup_{s\in(0,T]\cap I}\frac{1-e^{-rs}}{s}-r\right)+ \tilde D\cdot\left(\sup_{s\in(0,T]\cap I} \frac{{\gamma_1}^se^{-rs}-1}{s}-\left(\ln\gamma_1-r\right)\right)=C_0.$$ This gives a uniform pointwise estimate for the nonnegative function $\chi_{\left\{e^{rs}g>P_s(g\vee 0)\right\}}\left(g-e^{-r{s}}P_{{s}}g\right)$ on $E^s$ from which the general case of the Lemma's estimate can be derived immediately. If there exists an $i_0$ as stipulated in the Lemma, then by Remark \ref{Etproperties}, $$\bigcap_{s\in (0,T]\cap I}E^s = E^h$$ which completes the proof.
\end{proof}

Similarly, one can prove the corresponding estimate for calls with dividends (the dividends being assumed to be encoded in the discount rate $r$ and the Markov chain $P_\cdot$):

\begin{lem} \label{callestimateonE} Let this time $g=\bar f-K$ and suppose there is a $\gamma_0>0$ such that $$P_t\bar f\geq {\gamma_0}^t\bar f$$ for all $t\in(0,T]\cap I$. Let us define $$\bar D:=\chi_{(0,e^r)}\left(\gamma_0\right)\sup_{\bigcup_{t\in(0,T]\cap I}E^t}\bar f +\chi_{(e^r,+\infty)}\left(\gamma_0\right)\inf_{\bigcup_{t\in(0,T]\cap I}E^t}\bar f\geq 0.$$
Then there is a constant $C_1\in \RR$ given by $$C_1:=K\left(\sup_{s\in(0,T]\cap I}\frac{1-e^{-rs}}{s}- r\right)+\bar D\cdot \left(\sup_{s\in(0,T]\cap I} \frac{{\gamma_0}^se^{-rs}-1}{s}-\ln\gamma_0+ r\right)$$ such that for all $s\in(0,T]\cap I$ and measurable $A$, \begin{eqnarray*}&& \left\|\chi_{\left\{g> e^{-rs}P_{s}(g\vee 0)\right\}}\left(g-e^{-r{s}}P_{{s}}(g\vee 0)\right) \right\|_{L^1\left(\bigcap_{s\in (0,T]\cap I}E^s\cap A\right)} \\ &\leq &\left\|\chi_{\left\{g> e^{-rs}P_{s}(g\vee 0)\right\}}\left(g-e^{-r{s}}P_{{s}}(g\vee 0)\right) \right\|_{L^1(E^s\cap A)} \\ &\leq& \lambda^d\left[\left\{e^{rs}g>P_s(g\vee 0)>P_s g\right\}\cap A\right] \cdot \left(\left(\ln\gamma_0-r\right)\bar D+rK+ C_1 \right)\cdot s .\end{eqnarray*}
If moreover, there exists an $i_0\in\{1,\dots,m\}$ such that for any $j\in\{1,\dots, d\}$, $\left(y_{i_0}^{(h)}\right)_j\geq 0 $ for all $j\in\{1,\dots, d\}$, then we even have for all measurable $A$ \begin{eqnarray*}&& \left\|\chi_{\left\{g> e^{-rs}P_{s}(g\vee 0)\right\}}\left(g-e^{-r{s}}P_{{s}}(g\vee 0)\right) \right\|_{L^1\left(E^h\cap A\right)}  \\ &\leq& \lambda^d\left[\left\{e^{rs}g>P_s(g\vee 0)>P_s g\right\}\cap A\right] \cdot \left(\left(\ln\gamma_0-r\right)\bar D+rK+ C_1 \right)\cdot s .\end{eqnarray*}
\end{lem}


Next one will endeavour to find estimates on the Lebesgue measure of the set $\left\{e^{rs}g>P_s(g\vee 0)>P_s g\right\}\cap A$ occurring on the right hand side of the previous Lemmas \ref{estimateonE} and \ref{callestimateonE}.

\begin{lem} \label{setestimate}Suppose either 
\begin{enumerate}
\item $g$ is monotonely decreasing, eg $g=K-\bar f$, and
\item there is an $i_0\in\{1,\dots,m\}$ such that $y_{i_1}^{(h)}\leq 0$ componentwise 
\end{enumerate}
or one has
\begin{enumerate}
\item $g$ is monotonely increasing, eg $g=\bar f- K$, and
\item there exists an $i_1\in\{1,\dots,m\}$ such that $y_{i_1}^{(h)}\geq 0$ componentwise.
\end{enumerate}
Then for all measurable $A\subseteq\RR^d$,
\begin{eqnarray*}&&   \lambda^d\left[\left\{e^{rs}g>P_s(g\vee 0)>P_s g\right\}\cap A\right]\\ &\leq& \lambda^d\left[\left\{P_s(g\vee 0)>P_s g>0\right\}\cap A\right] + \lambda^d\left[\left\{P_s(g\vee 0)>0\right\}\cap\left\{P_s g\leq 0\right\}\cap A\right]\end{eqnarray*} 
\end{lem}

\begin{proof} We shall establish an upper bound for the set $\left\{e^{rs}g>P_s(g\vee 0)>P_sg\right\}$. Our assumptions entail that $$\left\{g>0\right\}\subseteq \left\{P_s(g\vee 0)>0\right\}$$ for all $s\in(0,T]\cap I$. From this we may, using $P_s(g\vee 0)\geq 0$, derive \begin{eqnarray*} && \left\{e^{rs}g>P_s(g\vee 0)>P_sg\right\}\cap \left\{P_sg\leq 0\right\} \\ &\subseteq & \{g>0\}\cap \left\{P_s(g\vee 0)>P_sg\right\}\cap \left\{P_sg\leq 0\right\}\\ &\subseteq& \left\{P_s(g\vee 0)>0\right\}\cap \left\{P_sg\leq 0\right\}\cap \left\{P_s(g\vee 0)>P_sg\right\} \\ &=& \left\{P_s(g\vee 0)>0\right\}\cap \left\{P_sg\leq 0\right\}\end{eqnarray*} This implies \begin{eqnarray*} && \left\{e^{rs}g>P_s(g\vee 0)>P_sg\right\}\\ &\subseteq& \left\{P_s(g\vee 0)>P_sg>0\right\}\cup \left(\left\{P_s(g\vee 0)>0\right\}\cap \left\{P_sg\leq 0\right\}\right). \end{eqnarray*}
\end{proof}

The estimate of the preceding Lemma will become relevant thanks to the following result (which in turn is based on the Corollary \ref{Bs/2-Bsequation} and Lemma \ref{estimateonE}).

\begin{lem} \label{chapter8a,Lemma1.9}Suppose $ y_i^{\left(h\right)}\geq 0$ componentwise for all $i\in\{1,\dots, m\}$. Then for all $s\in(0,T]\cap (2\cdot I)$ and $f\geq 0$, \begin{eqnarray*} &&\left\|\left(B_{\frac{s}{2}}\right)^{\circ 2}f-B_sf\right\|_{L^1\left(E^h \cap \bigcap_{k=1}^{m^\frac{s}{2h}} \left(A -y_k^{\left(\frac{s}{2}\right)}\right)\right)} \\ & \leq & \left\|\left(B_{\frac{s}{2}}\right)^{\circ 2}f-B_sf\right\|_{L^1\left(E^\frac{s}{2}\cap \bigcap_{k=1}^{m^\frac{s}{2h}} \left(A-y_k^{\left(\frac{s}{2}\right)}\right)\right)}\end{eqnarray*}

Assume furthermore $g=K-\bar f$ and there exists a real number $\gamma_1>0$ such that $P_h\bar f\leq{\gamma_1}^h\bar f$. Then for all $f\geq 0$, \begin{eqnarray*} &&\left\|\left(B_{\frac{s}{2}}\right)^{\circ 2}f-B_sf\right\|_{L^1\left(E^h \cap \bigcap_{k=1}^{m^\frac{s}{2h}} \left(A -y_k^{\left(\frac{s}{2}\right)}\right)\right)} \\ & \leq & \left\|\left(B_{\frac{s}{2}}\right)^{\circ 2}f-B_sf\right\|_{L^1\left(E^\frac{s}{2}\cap \bigcap_{k=1}^{m^\frac{s}{2h}} \left(A-y_k^{\left(\frac{s}{2}\right)}\right)\right)} \\ &\leq& \lambda^d\left[\left\{e^{rs/2}g>P_{s/2}(g\vee 0)>P_{s/2} g\right\}\cap A\right]\cdot e^{-r\frac{s}{2}}\\ &&\cdot \left(\left(\ln\gamma_1-r\right)\tilde D +rK+ C_0  \right)\cdot \frac{s}{2} \end{eqnarray*} with $C_0$ and $\tilde D$ as before.

\end{lem}
\begin{proof} Consider $t\in I$. Due to our assumption of $\left(\min_{i\in\{1,\dots,m\}} \left(y_i^{\left(h\right)}\right)_j\right)_{j\in\{1,\dots,d\}}\geq 0$ componentwise (which according to the notational convention introduced at the outset, can be written $\min_{i\in\{1,\dots,m\}} y_i^{\left(h\right)}\geq 0$), one has $$\min_k y_k^{\left(t\right)}=\frac{t}{h}\min_{i} y_i^{\left(h\right)}\geq 0.$$ Since the set $E^t$ is north-east connected, this yields ${E^t-y_k^{\left(t\right)}}\supseteq {E^t}$ for all $k\in\left\{1,\dots,m^{\frac{t}{h}}\right\}$, which in turn -- via $\chi_{E^t}\left(\cdot{+y_k^{\left(t\right)}}\right)=\chi_{E^t-y_k^{\left(t\right)}}\geq \chi_{E^t}$ for all $k\in\left\{1,\dots, m^\frac{t}{h}\right\}$ -- gives 
\begin{eqnarray*} P_{t}\left(\chi_{E^t}f\right)&=& \sum_{k=1}^{m^\frac{t}{h}}\alpha^{(t)}_k\chi_{E^t} \left(\cdot+y_k^{(t)}\right)f\left(\cdot+y_k^{(t)}\right)\\ &\geq& \sum_{k=1}^{m^\frac{t}{h}}\alpha^{(t)}_k\chi_{E^ t} \left(\cdot\right)f\left(\cdot+y_k^{(t)}\right)\\ &=& \chi_{E^t}{P_{t}}f\end{eqnarray*} for all $f\geq 0$. The same holds of course when replacing $f$ by $f\chi_A$. Moreover, treating the cases of $\min_k\chi_{A}\left(\cdot+y_k^{\left(t\right)}\right)=0$ and $\min_k\chi_{A}\left(\cdot+y_k^{\left(t\right)}\right)=1$ separately, we also see that $$\left(\min_k\chi_{A}\left(\cdot+y_k^{\left(t\right)}\right)\right){P_t}f\leq{P_t}\left(\chi_Af\right)$$ for all $f\geq 0$. Summarising these observations, we obtain
\begin{eqnarray*}P_{t}\left(\chi_{E^t\cap A}f\right) = P_{t}\left(\chi_{E^t}\cdot\chi_A f\right)&\geq& \chi_{E^t}\cdot{P_t}\left(\chi_Af\right) \\&\geq& \chi_{E^t}\cdot\left(\min_k\chi_{A}\left(\cdot+y_k^{\left(t\right)}\right)\right){P_t}f\\&\geq & \chi_{E^t}\chi_{\bigcap_k\left(A-y_k^{\left(t\right)}\right)} {P_t}f\end{eqnarray*} for all $f\geq 0$. Therefore -- using in addition the translation-invariance of $P_{t}$ and $\lambda^d$ (which makes $P_t$ a map that preserves the $L^1\left(\lambda^d\right)$-norm of nonnegative measurable functions) -- we deduce that for all measurable $f\geq 0$, 
\begin{eqnarray*} \left\|{P_t}f\right\|_{L^1\left({E^t}\cap \bigcap_{k=1}^{m^\frac{t}{h}} \left(A-y_k^{\left(t\right)}\right)\right)}&\leq& \left\| P_t\left(f\chi_{E^t\cap A}\right) \right\|_{L^1\left(\RR^d\right)}\\ &=&\left\| f\chi_{E^t\cap A} \right\|_{L^1\left(\RR^d\right)} =\left\| f \right\|_{L^1\left(E^t\cap A\right)}.\end{eqnarray*} 
From this, using Corollary \ref{Bs/2-Bsequation}, we obtain
\begin{eqnarray*} &&\left\|\left(B_{\frac{s}{2}}\right)^{\circ 2}f-B_sf\right\|_{L^1\left(E^\frac{s}{2}\cap \bigcap_k\left(A-y_k^{\left(\frac{s}{2}\right)}\right)\right)}\\ &\leq& \left\|e^{-r\frac{s}{2}}P_{\frac{s}{2}}\left(\chi_{\left\{g>e^{-r\frac{s}{2}}P_{\frac{s}{2}}(g\vee 0)\right\}}\cdot\left(g-e^{-r\frac{s}{2}}P_{\frac{s}{2}}(g\vee 0)\right) \right)\right\|_{L^1\left(E^\frac{s}{2}\cap \bigcap_k\left(A-y_k^{\left(\frac{s}{2}\right)}\right) \right)} \\ &\leq& e^{-r\frac{s}{2}} \left\| \chi_{\left\{g>e^{-r\frac{s}{2}}P_{\frac{s}{2}}(g\vee 0)\right\}}\cdot\left(g-e^{-r\frac{s}{2}}P_{\frac{s}{2}}(g\vee 0)\right) \right\|_{L^1\left(E^\frac{s}{2}\cap A\right)} \end{eqnarray*}
This is enough to prove the general part of the Lemma (which holds for arbitrary $g$) -- in the situation of $g=K-\bar f$ with at most exponentially increasing $t\mapsto\frac{ P_t\bar f}{\bar f}$ (on $\{\bar f>0\}$), one can take advantage of Lemma \ref{estimateonE} to complete the proof of Lemma \ref{chapter8a,Lemma1.9}.
\end{proof}

In the situation of $d=1$, $g=K-\bar f$ and $\bar f=\exp$ (vanilla put), Lemma \ref{chapter8a,Lemma1.9} can be readily combined with Lemma \ref{chapter8a,1dvanilla} to prove Theorem \ref{estimateBs/2-BsonE}, an estimate on the $L^1\left(E^h\right)$-norm of the difference $\left(B_{\frac{s}{2}}\right)^{\circ 2}f - B_sf$ that is quadratic in $s$.

\begin{rem} Let the translation-invariant Markov semigroup $P$ be derived from a cubature formula for the Gaussian measure with points $\{z_1,\dots,z_m\}$ in such a way that a geometric Brownian motion with logarithmic drift $\mu=\left(r-\frac{{\sigma_k}^2}{2}\right)_{k\in\{1,\dots,d\}}$ ($r>0$ and $\sigma\in{\RR_+}^d$ being the interest rate of the price process and the volatility vector, respectively) shall be approximated, that is to say $$\forall i\in\{1,\dots,m\}\forall k\in\{1,\dots,d\}\quad \left(y_i^{(h)}\right)_k=\mu_k h+\sigma_k h^{\frac{1}{2}}\left(z_i\right)_k.$$ Then the assumption that all the $y_i^{(h)}$ be componentwise nonnegative for $i\in\{1,\dots,m\}$ reads $$\min_{i} y_i^{\left(h\right)}=\mu h+\left(\min_i\left(z_i\right)_k\cdot \sigma_k h^\frac{1}{2}\right)_{k\in\{1,\dots,d\}}\geq 0$$ and therefore simply means that $\mu h$ is componentwise at least as big or even bigger than $-\left(\min_i\left(z_i\right)_k\cdot \sigma_k h^\frac{1}{2}\right)_{k\in\{1,\dots,d\}}$ which, needless to say, equals \linebreak $\left(\max_{i}\left(z_i\right)_k\cdot \sigma_k h^\frac{1}{2}\right)_{k\in\{1,\dots,d\}}$ in case of an axis-symmetric cubature formula, eg a cubature formula for the normal Gaussian measure. This assumption is tantamount to $$\forall k\in\{1,\dots,d\}\quad {\sigma_k}^2-2h^{-\frac{1}{2}}\min_i\left(z_i\right)_k\cdot\sigma_k -2{r}\leq 0, $$ which means (because of $h,r>0$ and $\sigma_k\ge 0$ for all $k$): \begin{eqnarray*}\forall k\in\{1,\dots,d\}\quad {\sigma_k}&\in& \RR_+\cap\left( \begin{array}{c} h^{-\frac{1}{2}}\min_i\left(z_i\right)_k-\sqrt{h^{-1}\cdot\left(\min_i\left(z_i\right)_k\right)^2 + 2{r}}, \\ h^{-\frac{1}{2}}\min_i\left(z_i\right)_k +\sqrt{h^{-1}\cdot\left(\min_i\left(z_i\right)_k\right)^2 + 2{r}}\end{array}\right)\\ &&= h^{-\frac{1}{2}}\cdot \left[0, \min_i\left(z_i\right)_k + \sqrt{\left(\min_i\left(z_i\right)_k\right)^2 + 2rh}\right]\end{eqnarray*} for all $k\in\{1,\dots,d\}$, entailing that $P$ models a basket of logarithmic asset prices whose volatilities are bounded above by the positive number \linebreak $h^{-\frac{1}{2}}\left(\min_i\left(z_i\right)_k +\sqrt{\left(\min_i\left(z_i\right)_k\right)^2 + 2rh}\right)$.

Later on, in Example \ref{cubature example}, we shall present an analogous reasoning under the assumption that $\mu$ is chosen directly from the condition that $e^{-(r-\delta)\cdot}P_\cdot\bar f$ -- wherein $\delta$ denotes the continuous dividend yield -- be a martingale.
\end{rem}

Now, owing to the pecularity that our investigations are only concerned with discrete translation-invariant Markov chains $\left(P_t\right)_{t\in I}$ (Markov chains which are derived from cubature formulae, for instance), we can use rather elementary inequalities to find upper bounds on the subsets of $\RR^d$ occurring in the estimates of Lemma \ref{estimateonE}.

As mentioned previously, we will start with the simple, nevertheless practically important, example of a one-dimensional American vanilla put:

\begin{lem} \label{chapter8a,1dvanilla}Suppose $d=1$ and $\bar f=\exp$, and let $g=K-\bar f$. Under these assumptions there exists a $\gamma_1>0$ such that $P_t\bar f={\gamma_1}^t\bar f$, and furthermore, one has for all $s\in I$, $$\left\{e^{rs}g>P_s(g\vee 0)>P_s g\right\}\subseteq\ln K+ \left( -\frac{1}{h}\cdot\max_{i\in\left\{1,\dots,m \right\}} y_i^{(h)}, 0\right)\cdot s.$$
\end{lem}
\begin{proof} The real number $\gamma_1$ is given by the relation $${\gamma_1}^h=\sum_{i=1}^m \alpha_i^{(h)}e^{y_i^{(h)}},$$ that is $$\gamma_1=e^\frac{\ln \left(\sum_{i=1}^m \alpha_i^{(h)}e^{y_i^{(h)}}\right)}{h}.$$
Next we observe that on the one hand by Remark \ref{Etproperties} \begin{eqnarray*}\left\{P_s(g\vee 0)>P_s g\right\}&=& \left\{\exists k\in\left\{1,\dots,m^\frac{s}{h}\right\}\quad g\left(\cdot+y_k^{(s)}\right)<0\right\} \\  &= &\left\{\min_{k\in\left\{1,\dots,m^\frac{s}{h}\right\}} g\left(\cdot+y_k^{(s)}\right)<0\right\} \\ &=& \left\{K-\max_{k\in\left\{1,\dots,m^\frac{s}{h}\right\}} \exp\left(\cdot+y_k^{(s)}\right)<0\right\}\\ &= &\left\{K- \exp\left(\cdot+\max_{k\in\left\{1,\dots,m^\frac{s}{h}\right\}}y_k^{(s)}\right)<0\right\}\\ &=& \left\{\ln K< \cdot+\max_{k\in\left\{1,\dots,m^\frac{s}{h}\right\}}y_k^{(s)}\right\} \\&=&\left(\ln K - \max_{k\in\left\{1,\dots,m^\frac{s}{h}\right\}}y_k^{(s)}, +\infty\right) \\&=&\left(\ln K - \frac{s}{h}\cdot\max_{i\in\left\{1,\dots,m \right\}}y_i^{(h)}, +\infty\right) \end{eqnarray*} and secondly 
$$\{g>0\}=\left\{K>\exp\right\}=\left(-\infty, \ln K\right),$$ thus \begin{eqnarray*}\left\{e^{rs}g>P_s(g\vee 0)>P_s g\right\}&\subseteq &\left(\ln K - \frac{s}{h}\cdot\max_{i\in\left\{1,\dots,m \right\}}y_i^{(h)}, \ln K\right)\\ &=& \ln K+ \left( -\frac{1}{h}\cdot\max_{i\in\left\{1,\dots,m \right\}}y_i^{(h)}, 0\right)\cdot s \end{eqnarray*}
\end{proof}

In light of Lemma \ref{chapter8a,Lemma1.9}, Lemma \ref{chapter8a,1dvanilla} can now finally be applied to prove the following result:

\begin{Th}\label{estimateBs/2-BsonE} Let $g=K-\bar f$ and suppose $d=1$ as well as $\bar f=\exp$. Under these assumptions there is a $\gamma_1>0$ such that $P_t\bar f={\gamma_1}^t\bar f$ for all $t\in I$. Assume, moreover, that $y_i^{\left(h\right)}\geq 0$ componentwise for all ${i\in\{1,\dots,m\}}$. Then there is a real constant $D$ such that for all $s\in(0,T]\cap (2\cdot I)$ and for all $f\geq g\vee 0$, $$\left\|\left(B_{\frac{s}{2}}\right)^{\circ 2}f - B_sf\right\|_{L^1\left(E^h\right)}\leq \left\|\left(B_{\frac{s}{2}}\right)^{\circ 2}f - B_sf\right\|_{L^1\left(E^\frac{s}{2}\right)}\leq\frac{D}{2}\cdot{s}^{2}.$$ We can compute $D$ explicitly from the constants $C_0$ and $\tilde D$ of Lemma \ref{estimateonE} as $$D=\left(\left(\ln\gamma_1-r\right)\tilde{D} +rK+ C_0  \right)\cdot \frac{\max_{i} y_i^{(h)}}{h}. $$
\end{Th}
\begin{proof} It is enough to combine Lemma \ref{chapter8a,Lemma1.9} for $A=\RR^d$ with Lemma \ref{chapter8a,1dvanilla}.
\end{proof}

\section{Bounds on higher order differences}

Now we shall proceed to establish convergence estimates for the sequence $\left(B_{T\cdot 2^{-n}}f\right)_{n\in\NN}$ in the $L^1(E^h\cap A)$-norm, for all measurable $f\geq g\vee 0$ and measurable $A\subseteq\RR^d$. Thus we have to prove bounds on the {\em higher order differences} $\left(B_{s\cdot 2^{-N}}\right)^{\circ \left(2^N\right)}f-\left(B_{s\cdot 2^{-M}}\right)^{\circ \left(2^M\right)}f$.

\begin{lem}\label{B_t-B_s_estimate} Suppose $d=1$ and $\bar f=\exp$. Under these assumptions there is a $\gamma_1$ such that $P_t\bar f={\gamma_1}^t\bar f$. Assume, moreover, that $y_i^{\left(h\right)}\geq 0$ componentwise for all ${i\in\{1,\dots,m\}}$. Under these assumptions there exists a real number $D>0$ (the same as in Theorem \ref{estimateBs/2-BsonE}) such that for all $k\in\NN_0$, $s\in(0,T]\cap \left(2^{k+1}\cdot I\right)$ and measurable $f\geq g\vee 0$, one has $$\left\|\left(B_{s\cdot 2^{-(k+1)}}\right)^{\circ\left(2^{k+1}\right)}f - \left(B_{s\cdot 2^{-k}}\right)^{\circ\left(2^k\right)}f\right\|_{L^1\left(\lambda^1\left[E^h\cap\cdot\right]\right)}\leq D\cdot {s}^2\cdot{2}^{-(k+1)}.$$ 
\end{lem}

The proof is contrived inductively, the base step being Theorem \ref{estimateBs/2-BsonE}, and the induction step being the first part of Lemma \ref{from_k=1_to_kinN}. However, the second and more general part of Lemma \ref{from_k=1_to_kinN} -- which we will need later on in this paper when we study options on multiple assets -- requires the following auxiliary result.

\begin{lem}\label{B_tcontractsL1Eh} Let $t\in I$, $A\subseteq\RR^d$ measurable, and assume also $y_i^{\left(h\right)}\geq 0$ componentwise for all ${i\in\{1,\dots,m\}}$. Then $\bigcap_{s\in (0,T]\cap I}E^s=E^h$ by Remark \ref{Etproperties}, and for all $f_1\geq f_0\geq g\vee 0$ and $p\in\{1,+\infty\}$, \begin{eqnarray*} && \left\|B_tf_1-B_tf_0\right\|_{L^p\left(\lambda^d\left[E^h\cap \bigcap_{k\in\left\{1,\dots, m^\frac{t}{h}\right\}}\left(A-y_k^{(t)}\right)\cap \cdot\right]\right)}\\ &\leq &e^{-rt}\left\|f_1-f_0\right\|_{L^p\left(\lambda^d\left[E^h\cap A\cap\cdot\right]\right)}\end{eqnarray*}
\end{lem}
\begin{proof} Consider a measurable set $A\subset \RR^d$ and measurable functions $f_0,f_1\geq g\vee 0$. Similarly to the proof of Lemma \ref{B_tcontractsifgeq g}, we observe that due to the monotonicity of $B_t$ and the fact that $B_tf\geq g\vee 0\geq g$ for all $f\geq 0$, \begin{eqnarray*}\left\{B_tf_1=g\right\}&=&\left\{B_tf_0\leq B_tf_1=g\right\} \\ &=& \left\{g\leq B_tf_0\leq B_tf_1=g\right\}\\ &=& \left\{B_tf_1=g\right\}\cap \left\{B_tf_0=g\right\}\\ &\subseteq& \left\{B_tf_1-B_tf_0=0\right\},\end{eqnarray*} that is $$\left\{B_tf_1-B_tf_0\neq 0\right\} \subseteq \left\{B_tf_1\neq g\right\} = \left\{B_tf_1>g\right\}$$ Combining this with the monotonicity of $P_t$ as well as the fact that $E^t+y_i^{(t)}\subseteq E^t$ for all $i$ (which in turn is a consequence of the north-east connectedness of $E^t$ -- cf Remark \ref{Etproperties} -- and the assumption that $y_i^{(t)}\geq 0$ componentwise for all $i$), we obtain
\begin{eqnarray}&& 0\leq\left(B_tf_1-B_tf_0\right)\chi_{E^h \cap \bigcap_k\left(A-y_k^{(t)}\right)}\nonumber  \\ &= & \left(B_tf_1-B_tf_0\right)\chi_{E^h \cap \bigcap_k\left(A-y_k^{(t)}\right)\cap \left\{B_tf_1>g\right\}}\nonumber  \\ &= & \chi_{E^h \cap  \bigcap_k\left(A-y_k^{(t)}\right)\cap\left\{B_tf_1>g\right\}}\left(e^{-rt}P_tf_1\vee g-e^{-rt}P_tf_0\vee g\right) \nonumber \\ &=& \chi_{E^h \cap  \bigcap_k\left(A-y_k^{(t)}\right)\cap\left\{e^{-rt}P_tf_1>g\right\}}\left(e^{-rt}P_tf_1-e^{-rt}P_tf_0\vee g\right) \nonumber \\ &\leq&\chi_{E^h \cap \bigcap_k\left(A-y_k^{(t)}\right)\cap \left\{e^{-rt}P_tf_1>g\right\}}\left(e^{-rt}P_tf_1-e^{-rt}P_tf_0\right) \nonumber \\ &\leq& e^{-rt}\chi_{E^h \cap \bigcap_k\left(A-y_k^{(t)}\right)}\left(P_tf_1-P_tf_0\right) \nonumber \\ &=&  e^{-rt}\sum_{i=1}^{m^{\frac{t}{h}}}\alpha_i^{(t)}\chi_{E^h \cap \bigcap_k\left(A-y_k^{(t)}\right)}\left(f_1\left(\cdot+y_i^{(t)}\right)-f_0\left(\cdot+y_i^{(t)}\right)\right) \nonumber \\ &=&  e^{-rt}\sum_{i=1}^{m^{\frac{t}{h}}}\alpha_i^{(t)}\left(\chi_{\left(E^h \cap \bigcap_k\left(A-y_k^{(t)}\right)\right)+y_i^{(t)}}\left(f_1-f_0\right)\right)\left(\cdot+y_i^{(t)}\right) \nonumber \\ &=&  e^{-rt}\sum_{i=1}^{m^{\frac{t}{h}}}\alpha_i^{(t)}\left(\chi_{\left(E^h+y_i^{(t)}\right) \cap \bigcap_k\left(A-y_k^{(t)}+y_i^{(t)}\right)}\left(f_1-f_0\right)\right)\left(\cdot+y_i^{(t)}\right) \nonumber\end{eqnarray}

Now, since $$\bigcap_{k\in\left\{1,\dots,m^\frac{t}{h}\right\}}\left(A-y_k^{(t)}+y_i^{(t)}\right) \subseteq A$$ for all $i\in\left\{1,\dots,m^\frac{t}{h}\right\}$ and $f_1-f_0\geq 0$, this means
\begin{eqnarray} && \left(B_tf_1-B_tf_0\right)\chi_{E^h \cap \bigcap_k\left(A-y_k^{(t)}\right)}\nonumber  \\ &\leq& e^{-rt}\sum_{i=1}^{m^{\frac{t}{h}}}\alpha_i^{(t)}\left(\chi_{\left(E^h+y_i^{(t)}\right) \cap A}\left(f_1-f_0\right)\right)\left(\cdot+y_i^{(t)}\right) \label{pointwiseBtdifference} 
\end{eqnarray} 

Combining this pointwise estimate with the translation-invariance of the Lebesgue measure yields
\begin{eqnarray}&& \left\|B_tf_1-B_tf_0\right\|_{L^1\left(\lambda^d\left[E^h \cap \bigcap_k\left(A-y_k^{(t)}\right)\cap\cdot\right]\right)}\nonumber  \\ &= &\int_{\RR^d} \left(B_tf_1-B_tf_0\right)\chi_{E^h \cap \bigcap_k\left(A-y_k^{(t)}\right)}d\lambda^d \nonumber \\ &\leq& e^{-rt}\sum_{i=1}^{m^{\frac{t}{h}}}\alpha_i^{(t)}\int_{\RR^d}\left(\chi_{\left(E^h+y_i^{(t)}\right) \cap A}\left(f_1-f_0\right)\right)\left(\cdot+y_i^{(t)}\right) d\lambda^d \nonumber \\ &=& e^{-rt}\sum_{i=1}^{m^{\frac{t}{h}}}\alpha_i^{(t)}\int_{\RR^d}\left(\chi_{\left(E^h+y_i^{(t)}\right) \cap A}\left(f_1-f_0\right)\right) d\lambda^d \nonumber \\ &\leq& e^{-rt}\sum_{i=1}^{m^{\frac{t}{h}}}\alpha_i^{(t)}\int_{E^h\cap A}\left(f_1-f_0\right) d\lambda^d\nonumber \\ &\leq& \int_{\RR^d} e^{-rt}\chi_{E^h \cap A}\left(f_1-f_0\right)  d\lambda^d\nonumber \\ &=& e^{-rt}\left\| f_1-f_0\right\|_{L^1\left(E^h \cap A\right)}\nonumber, \end{eqnarray}
where we have used the inclusion $E^s+y_k^{(t)}\subseteq E^s$ which -- owing to the north-east connectedness of the sets $E^s$ and our assumption $\min_{i\in\{1,\dots,m\}}y^{(h)}\geq 0$ -- holds for arbitrary $k\in\left\{1,\dots,m^\frac{t}{h}\right\}$ and $s,t\in I$ as well as the assumption $f_1-f_0\geq0$.
Similarly, based on \eqref{pointwiseBtdifference}, the translation-invariance and the sub-linearity of the $\esssup_{\RR^d}$-norm imply for measurable $f_1\geq f_0\geq g\vee 0$
\begin{eqnarray}&& \left\|B_tf_1-B_tf_0\right\|_{L^\infty\left(\lambda^d\left[E^h \cap \bigcap_k\left(A-y_k^{(t)}\right)\cap\cdot\right]\right)}\nonumber  \\ &= &\esssup_{\RR^d} \left[ \left(B_tf_1-B_tf_0\right)\chi_{E^h \cap \bigcap_k\left(A-y_k^{(t)}\right)}\right] \nonumber  \\ &\leq& \esssup_{\RR^d}\left[ e^{-rt}\sum_{i=1}^{m^{\frac{t}{h}}}\alpha_i^{(t)}\left(\chi_{\left(E^h+y_i^{(t)}\right) \cap A}\left(f_1-f_0\right)\right)\left(\cdot+y_i^{(t)}\right) \right]\nonumber \\ &\leq& e^{-rt}\sum_{i=1}^{m^{\frac{t}{h}}}\alpha_i^{(t)}\esssup_{\RR^d}\left[\left(\chi_{\left(E^h+y_i^{(t)}\right) \cap A} \left(f_1-f_0\right)\right)\left(\cdot+y_i^{(t)}\right) \right]\nonumber \\ &\leq& e^{-rt}\sum_{i=1}^{m^{\frac{t}{h}}}\alpha_i^{(t)}\esssup_{\RR^d}\left[\chi_{\left(E^h+y_i^{(t)}\right) \cap A} \left(f_1-f_0\right)\right] \nonumber \\ &\leq& e^{-rt}\sum_{i=1}^{m^{\frac{t}{h}}}\alpha_i^{(t)}\esssup_{\RR^d}\left[\chi_{E^h \cap A} \left(f_1-f_0\right)\right]\nonumber \\ &=& e^{-rt}\sum_{i=1}^{m^{\frac{t}{h}}}\alpha_i^{(t)}\esssup_{E^h \cap A} \left(f_1-f_0\right) \nonumber \\ &=& e^{-rt}\esssup_{E^h \cap A} \left(f_1-f_0\right) \nonumber  \end{eqnarray}
where again one has exploited the inclusion $E^s+y_k^{(t)}\subseteq E^s$ that holds for any $k\in\left\{1,\dots,m^\frac{t}{h}\right\}$ and $s,t\in I$.

\end{proof}

Lemma \ref{B_tcontractsL1Eh} plays a crucial r\^ole in the proof of the following Lemma \ref{from_k=1_to_kinN}, which forms a bridge between bounds on first order differences that are quadratic in the mesh of $I$ and bounds on higher order differences that also are quadratic in the mesh of $I$.

\begin{lem}\label{from_k=1_to_kinN} Let $T\in I$ and $p\in\{1,+\infty\}$. Consider a real number $D'>0$ and a measurable set $C\subseteq\RR^d$. Suppose one has an estimate of the kind \begin{eqnarray*}&&\forall f\geq g\vee 0\forall s\in (2\cdot I)\cap(0,T)\\ && \left\|\left(B_{\frac{s}{2}}\right)^{\circ 2}f - B_sf\right\|_{L^p\left(\lambda^d\left[C \cap\cdot\right]\right)}\leq \frac{D'}{2} \cdot{s}^{2}.\end{eqnarray*} Assume, moreover, $y_i^{\left(h\right)}\geq 0$ componentwise for all ${i\in\{1,\dots,m\}}$ (which by Remark \ref{Etproperties} also entails $\bigcap_{t\in (0,T]\cap I}E^t=E^h$). Then we get for all measurable $f\geq g\vee 0$ and for all $k\in\NN_0$, $s>0$  such that $s\in(0,T)\cap \left(2^{k+1}\cdot I\right)$, the estimate \begin{eqnarray*}&&\left\|\left(B_{s\cdot 2^{-(k+1)}}\right)^{\circ\left(2^{k+1}\right)}f - \left(B_{s\cdot 2^{-k}}\right)^{\circ\left(2^k\right)}f\right\|_{L^p\left(\lambda^d\left[C \cap\cdot\right]\right)}\\&\leq& D'\cdot {s}^2\cdot{2}^{-(k+1)}.\end{eqnarray*} Furthermore, if one assumes in addition $0\in\left\{ y_i^{(h)}\ : \ i\in\{1,\dots,m\}\right\}$, then one has a related implication for $L^p\left(E^h\cap \bigcap_{i\in\left\{1,\dots,m^{\frac{s}{2}}\right\}}\left(A-y_i^{\left(\frac{s}{2}\right)}\right)\right)$ instead of $L^p\left(C \right)$ for all measurable $A\subset\RR^d$: If under these assumptions the assertion \begin{eqnarray*}&&\forall f\geq g\vee 0\forall s\in (2\cdot I)\cap(0,T)\\ &&\left\|\left(B_{\frac{s}{2}}\right)^{\circ 2}f - B_sf\right\|_{L^p\left(E^h\cap \bigcap_{i}\left(A-y_i^{\left(s\right)}\right)\right)} \leq \frac{D'}{2} \cdot{s}^{2}\end{eqnarray*} holds, then the estimate \begin{eqnarray*}&& \forall f\geq g\vee 0\\ && \left\|\left(B_{s\cdot 2^{-(k+1)}}\right)^{\circ\left(2^{k+1}\right)}f - \left(B_{s\cdot 2^{-k}}\right)^{\circ\left(2^k\right)}f\right\|_{L^p\left(E^h \cap \bigcap_{i}\left(A-y_i^{\left(s\right)}\right)\right)}\\ &\leq &D'\cdot {s}^2\cdot{2}^{-(k+1)}\end{eqnarray*} is also valid for all $k\in\NN_0$ and $s>0$ such that $s\in(0,T)\cap \left(2^{k+1}\cdot I\right)$.

\end{lem}

\begin{proof} For both parts of the Lemma, we will conduct an induction in $k\in\NN_0$, the initial (or base) step being tautological each time. We have for all $s\in(0,T)\cap \left(2^{k+1}\cdot I\right)$ and $f\geq g\vee 0$ the estimate \begin{eqnarray*} &&\left(B_{s\cdot 2^{-(k+1)}}\right)^{\circ\left(2^{k+1}\right)}f - \left(B_{s\cdot 2^{-k}}\right)^{\circ2^k}f\\ &=& \left(B_{s\cdot 2^{-(k+1)}}\right)^{\circ\left(2^{k}\right)}\circ\left(B_{s\cdot 2^{-(k+1)}}\right)^{\circ\left(2^{k}\right)}f \\ && -\left(B_{s\cdot 2^{-(k+1)}}\right)^{\circ\left(2^{k}\right)}\circ\left(B_{s\cdot 2^{-k}}\right)^{\circ\left(2^{k-1}\right)}f \\ &&  + \left(B_{s\cdot 2^{-(k+1)}}\right)^{\circ\left(2^{k}\right)}\circ\left(B_{s\cdot 2^{-k}}\right)^{\circ\left(2^{k-1}\right)}f \\ && - \left(B_{s\cdot 2^{-k}}\right)^{\circ\left(2^{k-1}\right)}\circ \left(B_{s\cdot 2^{-k}}\right)^{\circ\left(2^{k-1}\right)}f \\ &=& \left(B_{s\cdot 2^{-(k+1)}}\right)^{\circ\left(2^{k}\right)}\circ\left(\left(B_{\frac{s}{2}\cdot 2^{-k}}\right)^{\circ\left(2^{k}\right)}-\left(B_{\frac{s}{2}\cdot 2^{-(k-1)}}\right)^{\circ\left(2^{k-1}\right)}\right)f \\ &&  + \left(\left(B_{\frac{s}{2}\cdot 2^{-k}}\right)^{\circ\left(2^{k}\right)} - \left(B_{ \frac{s}{2} \cdot 2^{-(k-1)}}\right)^{\circ\left(2^{k-1}\right)}\right)\circ \left(B_{s\cdot 2^{-k}}\right)^{\circ\left(2^{k-1}\right)}f \end{eqnarray*} which plays a crucial part in both the first and the second part of the Lemma.
For, we can first of all note that the induction hypothesis in the situation of the first part of the Lemma reads \begin{eqnarray}&&\forall f\geq g\vee 0\forall t\in (2^k\cdot I)\cap(0,T)\nonumber \\ &&  \left\| \left(B_{t\cdot 2^{-k}}\right)^{\circ\left(2^{k}\right)}f - \left(B_{ t\cdot 2^{-(k-1)}}\right)^{\circ\left(2^{k-1}\right)}f \right\|_{L^p\left(C\right)}\nonumber\\ &\leq& D'\cdot t^2\cdot 2^{-k}.\label{chapter8a,inductionhypothesis1}
\end{eqnarray} 
And if one now applies this induction hypothesis (\ref{chapter8a,inductionhypothesis1}) for $t=\frac{s}{2}$ (recalling that by assumption $s\in 2^{k+1}\cdot I$, thus $\frac{s}{2}\in 2^k\cdot I$) to the previous two equations and uses Lemma \ref{B_tcontractsL1Eh}, then one gets by the triangle inequality for the ${L^p\left(C\right)}$-norm,\begin{eqnarray*} && \left\|\left(B_{s\cdot 2^{-(k+1)}}\right)^{\circ\left(2^{k+1}\right)}f - \left(B_{s\cdot 2^{-k}}\right)^{\circ2^k}f\right\|_{L^p\left(C\right)} \\ &\leq & \left\|\left(\left(B_{\frac{s}{2}\cdot 2^{-k}}\right)^{\circ\left(2^{k}\right)} - \left(B_{ \frac{s}{2} \cdot 2^{-(k-1)}}\right)^{\circ\left(2^{k-1}\right)}\right)\circ \left(B_{s\cdot 2^{-k}}\right)^{\circ\left(2^{k-1}\right)}f \right\|_{L^p\left(C\right)}\\ && + \left\| \left(B_{s\cdot 2^{-(k+1)}}\right)^{\circ\left(2^{k}\right)}\circ\left(\left(B_{\frac{s}{2}\cdot 2^{-k}}\right)^{\circ\left(2^{k}\right)}-\left(B_{\frac{s}{2}\cdot 2^{-(k-1)}}\right)^{\circ\left(2^{k-1}\right)}\right)f \right\|_{L^p\left(C\right)}\\ &\leq & \left\|\left(\left(B_{\frac{s}{2}\cdot 2^{-k}}\right)^{\circ\left(2^{k}\right)} - \left(B_{ \frac{s}{2} \cdot 2^{-(k-1)}}\right)^{\circ\left(2^{k-1}\right)}\right)\circ \left(B_{s\cdot 2^{-k}}\right)^{\circ\left(2^{k-1}\right)}f \right\|_{L^p\left(C\right)} \\ && + \left\| \left(\left(B_{\frac{s}{2}\cdot 2^{-k}}\right)^{\circ\left(2^{k}\right)}-\left(B_{\frac{s}{2}\cdot 2^{-(k-1)}}\right)^{\circ\left(2^{k-1}\right)}\right)f \right\|_{L^p\left(C\right)} \\ &\leq & D'\cdot\frac{s^2}{4}\cdot 2^{-k}+D'\cdot\frac{s^2}{4}\cdot 2^{-k} = D'\cdot{s}^2\cdot 2^{-(k+1)}.\end{eqnarray*} In order to be entitled to apply Lemma \ref{B_tcontractsL1Eh} in this situation we have successively used the fact that $$\forall t\in I \forall \ell\geq g\vee 0 \quad B_t\ell\geq g\vee 0.$$ This completes the induction step for the first part of the Lemma. 

Turning to the proof of the second assertion in the Lemma (where in addition to $\min_{i\in\{1,\dots,d\}}y_i^{(h)}\geq 0$, $0\in\left\{ y_i^{(h)}\ : \ i\in\{1,\dots,m\}\right\}$ is assumed), we remark that \begin{eqnarray}&&\forall t\in I\forall k\leq\frac{\ln t-\ln h}{\ln 2} \nonumber \\ A(t)&:=& \bigcap_{i\in\left\{1,\dots,m^{\frac{t}{h}}\right\}} \left(A-y_i^{\left({t}\right)}\right)\nonumber \\&=& \bigcap_{\ell,k\in\left\{1,\dots,m^{\frac{t}{h}}\right\}} \left(A-y_k^{\left(\frac{t}{2}\right)}-y_\ell^{\left(\frac{t}{2}\right)}\right)\nonumber \\ &=& \left(A\left(\frac{t}{2}\right)\right)\left(\frac{t}{2}\right)\\&=& \label{A(t)properties} \bigcap_{i_1,\dots,i_{2^k}\in\left\{1,\dots,m^{\frac{t}{h}\cdot2^{-k-1}}\right\}} \left(A-y_{i_1}^{\left( t \cdot 2^{-(k+1)}\right)}-\cdots-y_{i_{2^k}}^{\left( {t} \cdot 2^{-(k+1)}\right)}\right).\end{eqnarray} In particular, if $0\in\left\{ y_i^{(h)}\ : \ i\in\{1,\dots,m\}\right\}$, $A(s)$ is decreasing in $s$: $$\forall s,t\in I\left(s\leq t\Rightarrow A(s)\supseteq A(t)\right).$$ Similarly to proof of the first part of the present Lemma, we deduce \begin{eqnarray} && \left\|\left(B_{s\cdot 2^{-(k+1)}}\right)^{\circ\left(2^{k+1}\right)}f - \left(B_{s\cdot 2^{-k}}\right)^{\circ2^k}f\right\|_{L^p\left(E^h \cap A\left({s}\right) \right)} \nonumber \\ &\leq & \left\|\left(B_{s\cdot 2^{-(k+1)}}\right)^{\circ\left(2^{k+1}\right)}f - \left(B_{s\cdot 2^{-k}}\right)^{\circ2^k}f\right\|_{L^p\left(E^h \cap A\left(s\right) \right)} \nonumber \\ &\leq & \left\|\left(\left(B_{\frac{s}{2}\cdot 2^{-k}}\right)^{\circ\left(2^{k}\right)} - \left(B_{ \frac{s}{2} \cdot 2^{-(k-1)}}\right)^{\circ\left(2^{k-1}\right)}\right)\circ \left(B_{s\cdot 2^{-k}}\right)^{\circ\left(2^{k-1}\right)}f \right\|_{L^p\left(E^h \cap A\left(s\right) \right)} \nonumber \\ \nonumber && + \left\| \left(B_{s\cdot 2^{-(k+1)}}\right)^{\circ\left(2^{k}\right)}\circ\left(\left(B_{\frac{s}{2}\cdot 2^{-k}}\right)^{\circ\left(2^{k}\right)}-\left(B_{\frac{s}{2}\cdot 2^{-(k-1)}}\right)^{\circ\left(2^{k-1}\right)}\right)f \right\|_{L^p\left(E^h \cap A\left(2^{k+1}\frac{s}{2}\right) \right)} \\ \label{1ast} &\leq & \left\|\left(\left(B_{\frac{s}{2}\cdot 2^{-k}}\right)^{\circ\left(2^{k}\right)} - \left(B_{ \frac{s}{2} \cdot 2^{-(k-1)}}\right)^{\circ\left(2^{k-1}\right)}\right)\circ \left(B_{s\cdot 2^{-k}}\right)^{\circ\left(2^{k-1}\right)}f \right\|_{L^p\left(E^h \cap A\left(s\right) \right)} \\\nonumber && + \left\| \left(B_{s\cdot 2^{-(k+1)}}\right)^{\circ\left(2^{k}\right)}\circ\left(\left(B_{\frac{s}{2}\cdot 2^{-k}}\right)^{\circ\left(2^{k}\right)}-\left(B_{\frac{s}{2}\cdot 2^{-(k-1)}}\right)^{\circ\left(2^{k-1}\right)}\right)f \right\|_{L^p\left(E^h \cap A\left(s\right)\right)} \end{eqnarray} from the triangle inequality.
But by a successive application of Lemma \ref{B_tcontractsL1Eh}, combined with the properties (\ref{A(t)properties}) of $A(\cdot)$, we have for all $f_1\geq f_0\geq g\vee 0$,
\begin{align}\nonumber & \left\|\left(B_{s\cdot 2^{-(k+1)}}\right)^{\circ\left(2^{k}\right)}\circ\left(f_1-f_0\right)\right\|_{L^p\left(E^h\cap A(s)\right)}\\ \nonumber \leq& \left\|\left(B_{s\cdot 2^{-(k+1)}}\right)^{\circ\left(2^{k}\right)}\circ\left(f_1-f_0\right)\right\|_{L^p\left(E^h\cap \bigcap_{\ell\in\left\{1,\dots,m^\frac{s}{2h}\right\}} \left( A\left(\frac{s}{2}\right) -y_\ell^{\left(\frac{s}{2}\right)}\right)\right)} \\ \nonumber  =&  \left\|\left(B_{s\cdot 2^{-(k+1)}}\right)^{\circ\left(2^{k}\right)}\circ\left(f_1-f_0\right)\right\|_{L^p\left(E^h\cap \bigcap_{i_1,\dots,i_{2^k}\in\left\{1,\dots,m^\frac{s}{2^{k+1}h}\right\}} \left(A\left(\frac{s}{2}\right)-y_{i_{1}}^{\left(\frac{s}{2^{k+1}}\right)}-\dots -y_{i_{2^k}}^{\left(\frac{s}{2^{k+1}}\right)}\right)\right)}\\ \nonumber \leq& \left\|\left(B_{s\cdot 2^{-(k+1)}}\right)^{\circ\left(2^{k}-1\right)}\circ\left(f_1-f_0\right)\right\|_{L^p\left(E^h\cap \bigcap_{i_1,\dots,i_{2^k}\in\left\{1,\dots,m^\frac{s}{2^{k+1}h}\right\}} \left(A\left(\frac{s}{2}\right)-y_{i_{1}}^{\left(\frac{s}{2^{k+1}}\right)}-\dots -y_{i_{2^k-1}}^{\left(\frac{s}{2^{k+1}}\right)}\right)\right)}\\ \nonumber  \leq& \\\nonumber  \vdots&\\ \leq& \left\|f_1-f_0\right\|_{L^p\left(E^h\cap A\left(\frac{s}{2}\right)\right)}\label{3ast}.\end{align} 

In light of the inclusion $A(s)\subseteq A\left(\frac{s}{2}\right)$, we finally obtain from combining estimates \eqref{1ast} and \eqref{3ast}
\begin{eqnarray*} && \left\|\left(B_{s\cdot 2^{-(k+1)}}\right)^{\circ\left(2^{k+1}\right)}f - \left(B_{s\cdot 2^{-k}}\right)^{\circ2^k}f\right\|_{L^p\left(E^h \cap A\left({s}\right) \right)}\\ &\leq & \left\|\left(\begin{array}{c}\left(B_{\frac{s}{2}\cdot 2^{-k}}\right)^{\circ 2^{k}} -\\ \left(B_{ \frac{s}{2} \cdot 2^{-(k-1)}}\right)^{\circ 2^{k-1}}\end{array}\right)\circ \left(B_{s\cdot 2^{-k}}\right)^{\circ 2^{k-1}}f \right\|_{L^p\left(E^h \cap A\left(s\right) \right)} \\ && + \left\| \left(\left(B_{\frac{s}{2}\cdot 2^{-k}}\right)^{\circ\left(2^{k}\right)}-\left(B_{\frac{s}{2}\cdot 2^{-(k-1)}}\right)^{\circ\left(2^{k-1}\right)}\right)f \right\|_{L^p\left(E^h \cap A\left(\frac{s}{2}\right) \right)}\\ &\leq& \left\|\left(\begin{array}{c}\left(B_{\frac{s}{2}\cdot 2^{-k}}\right)^{\circ 2^{k}} -\\ \left(B_{ \frac{s}{2} \cdot 2^{-(k-1)}}\right)^{\circ 2^{k-1}}\end{array}\right)\circ \left(B_{s\cdot 2^{-k}}\right)^{\circ 2^{k-1}}f \right\|_{L^p\left(E^h \cap A\left(\frac{s}{2}\right) \right)} \\ && + \left\| \left(\left(B_{\frac{s}{2}\cdot 2^{-k}}\right)^{\circ\left(2^{k}\right)}-\left(B_{\frac{s}{2}\cdot 2^{-(k-1)}}\right)^{\circ\left(2^{k-1}\right)}\right)f \right\|_{L^p\left(E^h \cap A\left(\frac{s}{2}\right) \right)} \\ &\leq & D'\cdot\frac{s^2}{4}\cdot 2^{-k}+D'\cdot\frac{s^2}{4}\cdot 2^{-k} = D'\cdot{s}^2\cdot 2^{-(k+1)},\end{eqnarray*} 
where in the last line we have taken advantage of the induction hypothesis \begin{eqnarray*}&&\forall k\in\NN_0\forall f\geq g\vee 0\forall t\in (2^k\cdot I)\cap(0,T) \\ &&\left(\begin{array}{c}\left\| \left(B_{t\cdot 2^{-k}}\right)^{\circ\left(2^{k}\right)}f - \left(B_{ t\cdot 2^{-(k-1)}}\right)^{\circ\left(2^{k-1}\right)}f \right\|_{L^p\left(E^h \cap \bigcap_{i}\left(A-y_i^{\left(2^{k} {t} \right)}\right)\right)}\\ \leq D'\cdot t^2\cdot 2^{-k} \end{array}\right)\end{eqnarray*} for the special case $t=\frac{s}{2}$
\end{proof}

The assumption of $0\in\left\{ y_1^{(h)},\dots, y_m^{(h)}\right\}\subseteq\RR^d$ while $\min_{i\in\{1,\cdots,m\}}y_i^{(h)}\geq 0$ componentwise corresponds, in the application that we have in mind, to the volatility attaining a certain critical value:

\begin{ex}\label{cubature example} As already hinted at when stating Remark \ref{treeispolynomial}, one could think of $P_\cdot$ as a discretisation of, say, the one-dimensional Black-Scholes model with constant discount rate $r$ and volatility $\sigma$, constructed via Gaussian cubature. Consider a cubature formula of a certain degree for the one-dimensional normal Gaussian measure $\nu_{0,1}=\frac{e^{-x^2/2}dx}{\sqrt{2\pi}}$ with cubature points $z_1,\dots,z_m$ and weights $\alpha^{(h)}_1,\dots,\alpha^{(h)}_m$ (this entails $\min_jz_j<0$) which will then give rise to a new Markov chain via $$\forall i\in\{1,\dots,m\} \quad y_i^{(h)} = \mu h+ z_i\sigma h^\frac{1}{2}$$ where $\mu$ is chosen such that the translation-invariant Markov chain with increments $y_1^{(h)},\dots, y_m^{(h)}$ and weights $\alpha^{(h)}_1,\dots,\alpha^{(h)}_m$ defines a process whose exponential becomes, after discounting inflation at rate $r$ and a continuous dividend yield at rate $\delta>0$, a martingale:
\begin{align*} \mu=& r-\delta -\frac{1}{h}\ln\left( \sum_{i=1}^m\alpha^{(h)}_i e^{z_i\sigma h^\frac{1}{2}} \right) =  r-\delta-\frac{1}{h}\left(\sigma h^\frac{1}{2}\min_j z_j +\ln\left( \sum_{i=1}^m\alpha^{(h)}_i e^{\left(z_i-\min_j z_j\right)\sigma h^\frac{1}{2}} \right)\right)\end{align*}
Then one will have \begin{align*} \min_j y_j^{(h)} =& \left(r-\delta\right)h - \ln\left(\sum_{i=1}^m\alpha^{(h)}_i { e^{\sigma h^\frac{1}{2}z_i}/e^{\sigma h^\frac{1}{2}\min_j z_j} } \right)\\ =& \left(r-\delta\right)h - \ln\left(\sum_{i=1}^m\alpha^{(h)}_i \exp\left({\sigma h^\frac{1}{2}\underbrace{\left(z_i-\min_jz_j\right)}_{\geq 0}} \right) \right)
\end{align*}
which may become zero and even positive if $\sigma h^\frac{1}{2}$ is sufficiently small -- for, $\delta\leq r$ will hold to exclude risk-less arbitrage opportunities. 
(If simply $\{z_1,z_2\}=\{\pm 1\}$, then this was a discrete model for a logarithmic asset price evolution that converges weakly to the Black-Scholes model with volatility $\sigma$, dividend yield $\delta$ and discount rate $r$ as $h\downarrow 0$.) An analogous remark can be made for higher dimensions: A non-ruinous dividend yield combined with a relatively small volatility for all assets in the baskets will make $\min_j y_j^{(h)}$ vanish or even rise above nought componentwise.
\end{ex}

With the first half of Lemma \ref{from_k=1_to_kinN}, we have completed the proof of Lemma \ref{B_t-B_s_estimate}. We shall now apply this result to finally get to a convergence bound for $\left(B_{T\cdot 2^{-n}}(g\vee 0)\right)_n$ -- which can be conceived of as a sequence of non-perpetual Bermudan option prices when successively halving the exercise mesh size.

\begin{lem}\label{from_k,k+1_to_M,N} Let $p\in[1,+\infty]$. Consider a real constant $D>0$ as well as a measurable set $C$ and a set $E$ of nonnegative measurable functions, and suppose one has an estimate of the kind \begin{eqnarray*}&&\forall k\in\NN_0\forall f\in E\forall s\in (2^{k+1}\cdot I)\cap(0,T)\\ &&  \left\|\left(B_{s\cdot 2^{-(k+1)}}\right)^{\circ\left(2^{k+1}\right)}f - \left(B_{s\cdot 2^{-k}}\right)^{\circ\left(2^k\right)}f\right\|_{L^p\left(\lambda^d\left[C\cap \cdot\right]\right)}\leq D \cdot {s}^2\cdot{2}^{-(k+1)} .\end{eqnarray*}
Then for all $N>M\in\NN$, $s\in(0,T)\cap \left(2^N\cdot I\right)$ and $f\in E$, the estimate \begin{eqnarray*}\left\|\left(B_{s\cdot 2^{-N}}\right)^{\circ\left(2^{N}\right)}f - \left(B_{s\cdot 2^{-M}}\right)^{\circ\left(2^M\right)}f\right\|_{L^p(C)}&\leq& D\cdot {s^2}\cdot{2}^{-M}\left(1- 2^{-(N-M-1)}\right) \\ &\leq& D\cdot {s^2}\cdot{2}^{-M}\end{eqnarray*}
holds.
\end{lem}
\begin{proof} With $M,N$, $s$, $f$ as in the statement of the Lemma, we obtain by the triangle inequality \begin{eqnarray*} && \left\|\left(B_{s\cdot 2^{-N}}\right)^{\circ\left(2^{N}\right)}f - \left(B_{s\cdot 2^{-M}}\right)^{\circ\left(2^M\right)}f\right\|_{L^p(C)}\\ &=&\left\|\sum_{k=M}^{N-1} \left(\left(B_{s\cdot 2^{-(k+1)}}\right)^{\circ\left(2^{k+1}\right)}f - \left(B_{s\cdot 2^{-k}}\right)^{\circ\left(2^k\right)}f\right)\right\|_{L^p(C)} \\ &\leq&\sum_{k=M}^{N-1}\left\| \left(\left(B_{s\cdot 2^{-(k+1)}}\right)^{\circ\left(2^{k+1}\right)}f - \left(B_{s\cdot 2^{-k}}\right)^{\circ\left(2^k\right)}f\right)\right\|_{L^p(C)} \\ &\leq & \sum_{k=M}^{N-1}D\cdot s^2\cdot 2^{-k-1}=D\cdot\frac{s^2}{2}\cdot\sum_{k=0}^{N-1-M}2^{-k}2^{-M}\\ &=& D\cdot\frac{s^2}{2}\cdot 2^{-M}\cdot \frac{1-2^{-(N-M-1)}}{1-2^{-1}}\leq D\cdot\frac{s^2}{2}\cdot 2^{-M}\cdot 2.\end{eqnarray*}

\end{proof}

Thus, if we combine this last Lemma \ref{from_k,k+1_to_M,N} with Lemma \ref{B_t-B_s_estimate} we arrive at 

\begin{Th} \label{cubatureBermudanconv} Suppose, as before, $d=1$ and $\bar f=\exp$, as well as $g=K-\bar f$. Under these assumptions there is a $\gamma_1$ such that $P_t\bar f={\gamma_1}^t\bar f$ for every $t\in I$, and let us suppose this $\gamma_1\in(0,e^r]$. Assume furthermore $y_i^{\left(h\right)}\geq 0$ componentwise for all ${i\in\{1,\dots,m\}}$. Then there exists a real number $D>0$ such that for all $N>M\in\NN$, $s\in(0,T]\cap \left(2^N\cdot I\right)$ and monotonely decreasing $f\geq g\vee 0$, one has \begin{eqnarray*}&& \left\|\left(B_{s\cdot 2^{-N}}\right)^{\circ\left(2^{N}\right)}f - \left(B_{s\cdot 2^{-M}}\right)^{\circ\left(2^M\right)}f\right\|_{L^1\left(E^h\right)}\\ &\leq& D\cdot {s}^2\cdot{2}^{-M}\left(1- 2^{-(N-M-1)}\right) \\ &\leq& D\cdot {s}^2\cdot{2}^{-M}.\end{eqnarray*} $D$ is the constant of Theorem \ref{estimateBs/2-BsonE}: $$D=\left(\left(\ln\gamma_1-r\right)\tilde{D} +rK+ C_0  \right)\cdot \frac{\max_{i} y_i^{(h)}}{h} .$$
\end{Th}

\section{Application to American $\max$-put options}

Analogously, we shall finally proceed to prove convergence of quadratic order in $s$ for $\bar f=\sum_{j=1}^d w_j\exp\left((\cdot)_j\right)$, where $w_1,\dots,w_d$ is a convex combination (the weights for a weighted average of the components/assets in a $d$-dimensional basket), as well as for the choices $\bar f=\min_{j\in\{1,\dots,d\}}\exp\left((\cdot)_j\right)$ and $\bar f=\max_{j\in\{1,\dots,d\}}\exp\left((\cdot)_j\right)$. However, this time, we shall employ different norms: $L^1\left(\lambda^d\left[E^h\cap A\cap\cdot\right]\right)$ for a compact subset $A\subset \RR^d$ such that $\lambda^d\left[\bigcap_{s\in (0,T]\cap I}E^s \cap A\right]\in (0,+\infty)$.

The first part of this endeavour will be to prove the applicability of Lemma \ref{estimateonE}.
\begin{lem} \label{barf=max_gamma1} If $\bar f=\max_{j\in\{1,\dots,d\}}\exp\left((\cdot)_j\right)$, then $$P_s\bar f\leq {\gamma_1}^s\bar f $$ where $$\gamma_1:=\left(\underbrace{ \left( \sum_{i=1}^{m^\frac{s}{h}} \alpha_i^{(s)} \max_{j\in\{1,\dots,d\}} e^{\left(y_i^{(s)}\right)_j} \right) }_{>0} \right) ^{\frac{1}{h}}.$$
\end{lem}

\begin{proof} We have for all $s\in I$ the estimate
\begin{eqnarray*}P_s\bar f&=&\sum_{i=1}^{m^\frac{s}{h}} \alpha_i^{(s)} \bar f\left(\cdot+y_i^{(s)}\right)\\ &=& \sum_{i=1}^{m^\frac{s}{h}} \alpha_i^{(s)} \max_{j\in\{1,\dots,d\}}\exp\left(\left(\cdot+y_i^{(s)}\right)_j\right)\\ &\leq& \sum_{i=1}^{m^\frac{s}{h}} \alpha_i^{(s)} \max_{\ell\in\{1,\dots,d\}} e^{\left(y_i^{(s)}\right)_\ell} \max_{j\in\{1,\dots,d\}} \exp\left(\left(\cdot\right)_j\right)\\ &=& \underbrace{ \sum_{i=1}^{m^\frac{s}{h}} \alpha_i^{(s)} \max_{\ell\in\{1,\dots,d\}} e^{\left(y_i^{(s)}\right)_\ell} }_{>0}\max_{j}\exp\left(\left(\cdot\right)_j\right), \end{eqnarray*} in particular the estimate holds for $s=h$. But this is 
to say $$P_h\bar f\leq {\gamma_1}^h\bar f ,$$ hence we have established the estimate in the Lemma for $s=h$. This readily suffices to prove the Lemma's assertion, as $\left(P_s\right)_{s\in I}$ is a Markov semigroup and by applying the Chapman-Komogorov equation inductively, $$\forall n\in\NN\quad P_{nh} \bar f = \underbrace{P_{h}\cdots P_{h}}_{n}\bar f\leq \underbrace{{\gamma_1}^h\cdots {\gamma_1}^h}_{n}\bar f={\gamma_1}^{hn}\bar f.$$
\end{proof}

Analogously, one can easily prove the existence of such a $\gamma_1>0$ as required by Lemma \ref{estimateonE} for $\bar f=\sum_{j=1}^d w_j\exp\left((\cdot)_j\right)$ for a convex combination $w_1,\dots,w_d$ and $\bar f=\min_{j\in\{1,\dots,d\}}\exp\left((\cdot)_j\right)$.

We next turn our attention to deriving upper bounds for the measures of the sets in the estimates of Lemma \ref{estimateonE} for the said example 
of $\bar f=\max_{j\in\{1,\dots,d\}}\exp\left((\cdot)_j\right)$. We continue to use the notation $I=h\NN_0$ and $$P_s: f\mapsto\sum_{i=1}^{m^\frac{s}{h}}\alpha_i^{(s)} f\left(\cdot+y_i^{(s)}\right),$$ where $\left(P_s\right)_{s\in I}=\left(P_{nh}\right)_{n\in \NN_0}=\left(\underbrace{P_{h}\cdots P_{h}}_{n}\right)_{n\in \NN_0}$ is the Markov chain generated by $P_h$.

\begin{lem} \label{barf=max set estimates} If $\bar f=\max_{j\in\{1,\dots,d\}}\exp\left((\cdot)_j\right)$, then for all $s\in I$, \begin{eqnarray*}&&\left\{\forall i\in\left\{1,\dots,m^\frac{s}{h}\right\}\quad g\left(\cdot+y_i^{(s)}\right)\leq 0\right\}\\ &\subset&\bigcup_{j=1}^d \left( \underbrace{ \RR\times\cdots\times\RR }_{j-1} \times \left[\ln K- \frac{s}{h} \max_{i\in\{1,\dots,m\}} \left(y_i^{(h)}\right)_j ,+\infty\right) \times \underbrace{ \RR\times\cdots\times\RR }_{d-j} \right)\end{eqnarray*} as well as \begin{eqnarray*}&& \left\{\forall i\in\left\{1,\dots,m^\frac{s}{h}\right\}\quad g\left(\cdot+y_i^{(s)}\right)\geq 0\right\}\\ &=& \bigotimes_{j=1}^d\left(-\infty,\ln K- \frac{s}{h} \max_{i\in\{1,\dots,m\}} \left(y_i^{(h)}\right)_j\right]\end{eqnarray*}
\end{lem}
\begin{proof} Let $s\in I$. Then 
\begin{eqnarray*}&& \left\{\forall i\in\{1,\dots,m\}\quad g\left(\cdot+y_i^{(s)}\right)\leq 0\right\} \\ &= &\left\{\max_{i\in\left\{1,\dots,m^\frac{s}{h}\right\}} g\left(\cdot+y_i^{(s)}\right)\leq 0\right\} \\ &=& \left\{K-\min_{i\in\left\{1,\dots,m^\frac{s}{h}\right\} } \max_{j\in\{1,\dots,d\}} \exp\left(\left(\cdot+y_i^{(s)}\right)_j\right)\leq 0\right\}\\  &\subseteq& \left\{K-\max_{i\in\left\{1,\dots,m^\frac{s}{h}\right\}} \max_{j\in\{1,\dots,d\}} \exp\left(\left(\cdot+y_i^{(s)}\right)_j\right)\leq 0\right\}\\  &=& \left\{K-\max_{j\in\{1,\dots,d\}}\max_{i\in\left\{1,\dots,m^\frac{s}{h}\right\} } \exp\left(\left(\cdot+y_i^{(s)}\right)_j\right)\leq 0\right\}\\  &=& \left\{K- \max_{j\in\{1,\dots,d\}}\exp\left(\left(\cdot\right)_j +\max_{i\in\left\{1,\dots,m^\frac{s}{h}\right\}} \left(y_i^{(s)}\right)_j\right) \leq 0\right\} \\ &=& \left\{\ln K\leq \max_{j\in\{1,\dots,d\}}\left( \left(\cdot\right)_j + \max_{i\in\left\{1,\dots,m^\frac{s}{h}\right\}} \left(y_i^{(s)}\right)_j \right)\right\} \\ &=& \bigcup_{j=1}^d\left\{\ln K-\max_{i\in\left\{1,\dots,m^\frac{s}{h}\right\}} \left(y_i^{(s)}\right)_j\leq (\cdot)_j \right\}\\ &=& \bigcup_{j=1}^d\left\{\ln K-\frac{s}{h}\max_{i\in\left\{1,\dots,m\right\}} \left(y_i^{(h)}\right)_j\leq (\cdot)_j \right\}, \end{eqnarray*} 
and also 
\begin{eqnarray*}&& \left\{\forall i\in\{1,\dots,m\}\quad g\left(\cdot+y_i^{(s)}\right)\geq 0\right\} \\ &= &\left\{\min_{i\in\left\{1,\dots,m^\frac{s}{h}\right\}} g\left(\cdot+y_i^{(s)}\right)\geq 0\right\} \\ &=& \left\{K-\max_{i\in\left\{1,\dots,m^\frac{s}{h}\right\}} \max_{j\in\{1,\dots,d\}} \exp\left(\left(\cdot+y_i^{(s)}\right)_j\right)\geq 0\right\} \\  &=& \left\{K- \max_{j\in\{1,\dots,d\}}\max_{i\in\left\{1,\dots,m^\frac{s}{h}\right\}} \exp\left(\left(\cdot+y_i^{(s)}\right)_j\right)\geq  0\right\}\\  &=& \left\{K- \max_{j\in\{1,\dots,d\}}\exp\left(\left(\cdot\right)_j +\max_{i\in\left\{1,\dots,m^\frac{s}{h}\right\}} \left(y_i^{(s)}\right)_j\right)\geq 0\right\} \\  &=& \left\{K\geq \exp\left(\max_{j\in\{1,\dots,d\}}\left(\left(\cdot\right)_j +\max_{i\in\left\{1,\dots,m^\frac{s}{h}\right\}} \left(y_i^{(s)}\right)_j\right)\right)\right\} \\ &=& \left\{\ln K \geq  \max_{j\in\{1,\dots,d\}} \left(\left(\cdot\right)_j + \max_{i\in\left\{1,\dots,m^\frac{s}{h}\right\}} \left(y_i^{(s)}\right)_j \right)\right\} \\ &=& \bigcap_{j=1}^d\left\{\ln K-\max_{i\in\left\{1,\dots,m^\frac{s}{h}\right\}} \left(y_i^{(s)}\right)_j\geq  (\cdot)_j \right\}\\ &=& \bigcap_{j=1}^d\left\{\ln K-\frac{s}{h}\max_{i\in\left\{1,\dots,m\right\}} \left(y_i^{(h)}\right)_j\geq (\cdot)_j \right\}. \end{eqnarray*} 
\end{proof}

\begin{cor}\label{Psgv0>Psg.max} If $\bar f=\max_{j\in\{1,\dots,d\}}\exp\left((\cdot)_j\right)$, then for all $s\in I$,
\begin{eqnarray*}&&\left\{P_s(g\vee 0)>P_sg\right\} \\ &=& \bigcup_{j=1}^d \left( \underbrace{ \RR\times\cdots\times\RR }_{j-1} \times \left(\ln K - \frac{s}{h} \max_{i\in\{1,\dots,m\}} \left(y_i^{(h)}\right)_j ,+\infty\right) \times \underbrace{ \RR\times\cdots\times\RR }_{d-j} \right)\end{eqnarray*}
\end{cor}
\begin{proof} Let $s\in I$. We simply remark that
\begin{eqnarray*}&&\left\{P_s(g\vee 0)>P_sg\right\}\\ &=&\left\{\exists i\in\left\{1,\dots,m^\frac{s}{h}\right\}\quad g\left(\cdot+y_i^{(s)}\right)<0\right\} \\ &=&\complement \left\{\forall i\in\left\{1,\dots,m^\frac{s}{h}\right\}\quad g\left(\cdot+y_i^{(s)}\right)\geq 0\right\} \end{eqnarray*} and apply Lemma \ref{barf=max set estimates}.
\end{proof}

These estimates lead to the following Corollary that will enable us -- under the assumption of $y_{i}^{(h)}\geq 0$ for all $i\in\{1,\dots,m\}$ and $0\in\left\{y_i{(h)} \ : \ i\in\{1,\dots,m\}\right\}$ (in order to be entitled to apply eg Lemma \ref{from_k=1_to_kinN}) -- to prove an $L^1$-convergence estimate (on a particular subset of $\RR^d$) for $\left(B_{T\cdot 2^{-n}}f\right)_{n\in\NN}$ for any measurable $f\geq g\vee 0$.

\begin{cor} Suppose $\bar f=\max_{j\in\{1,\dots,d\}} \exp\left((\cdot)_j\right)$ and consider any compact set $B\subset \RR^d$. Then there is some $R>0$ such that $B-(\ln K)_{j=1}^d\subset\left[-{R},R\right]^d$. What is more, \begin{eqnarray*}&&\lambda^d\left[\left\{e^{rs}g>P_s(g\vee 0)>P_sg\right\}\cap B\right] \\&\leq& s\cdot R^{d-1}\frac{1}{h} \sum_{j=1}^d\left( \max_{i\in\{1,\dots,m\}} \left(y_i^{(h)}\right)_j \vee 0\right)\end{eqnarray*} for all $s\in I$.
\end{cor}
\begin{proof} Let $s\in I$. Since $$\left\{e^{rs}g>P_s(g\vee 0)>P_sg\right\}\subseteq \left\{P_s(g\vee 0)>P_sg\right\}\cap \{g>0\}$$ by the monotonicity of $P_s$, we only have to observe that \begin{eqnarray*} \{g>0\} &=& \left\{K>\max_{j\in\{1,\dots,d\}} \exp\left((\cdot)_j\right)\right\} \\ &=&\bigcap_{j=1}^d\left\{K>\exp\left((\cdot)_j\right)\right\}\\ &=&\bigotimes_{j=1}^d \left(-\infty,\ln K\right)\end{eqnarray*} to arrive -- after taking advantage of the preceding Corollary \ref{Psgv0>Psg.max} -- at \begin{eqnarray*}&& \left\{e^{rs}g>P_s(g\vee 0)>P_sg\right\}\\ &\subseteq & \bigcup_{j=1}^d \left( \begin{array}{c}\underbrace{ \left(-\infty,\ln K\right)\times\cdots\times\left(-\infty,\ln K\right) }_{j-1} \\\times \left(\ln K - \frac{s}{h} \max_{i\in\{1,\dots,m\}} \left(y_i^{(h)}\right)_j ,\ln K\right)\\ \times \underbrace{ \left(-\infty,\ln K\right)\times\cdots\times\left(-\infty,\ln K\right) }_{d-j}\end{array} \right) \\ &=& \bigcup_{j=1}^d \left( \underbrace{ \RR_{<0}\times\cdots\times\RR_{<0} }_{j-1} \times \left(-\frac{1}{h} \max_{i\in\{1,\dots,m\}} \left(y_i^{(h)}\right)_j,0 \right)\cdot s\times \underbrace{ \RR_{<0}\times\cdots\times\RR_{<0} }_{d-j} \right) \\ && +(\ln K)_{j=1}^d.\end{eqnarray*}

However, $R>0$ has been chosen such that $B-(\ln K)_{j=1}^d\subseteq\left[-{R},R\right]^d$. Thus \begin{eqnarray*}&&\left\{e^{rs}g>P_s(g\vee 0)>P_sg\right\}\cap B \\ &\subseteq& \bigcup_{j=1}^d \left( [-R,0)^{j-1} \times \left( -\frac{1}{h} \max_{i\in\{1,\dots,m\}} \left(y_i^{(h)}\right)_j,0 \right)\cdot s\times \RR^{d-j} \right) \\ && +(\ln K)_{j=1}^d,\end{eqnarray*} and from this inclusion we may deduce the estimate given in the Lemma.
\end{proof}

The inequality we have just derived implies that the \linebreak $\lambda^d\left[E^h \cap \bigcap_{\ell}\left(B-y_\ell^{\left(s\right)}\right)\cap\cdot\right]$-volume of the set occurring in Lemma \ref{estimateonE} is of order $s$ for any compact $B$ and for $\bar f=\max_{j\in\{1,\dots,d\}} \exp\left((\cdot)_j\right)$. Hence again by Lemma \ref{estimateonE} (which is applicable because of Lemma \ref{barf=max_gamma1}) we obtain that the difference $B_{s/2}f-B_sf$ is of order $s^2$ (this time, however in the $L^1\left(E^h \cap \bigcap_{\ell}\left(B-y_\ell^{(s)}\right)\right)$-norm). This estimate on the norm of $B_{s/2}f-B_sf$ leads, via Lemmas \ref{from_k,k+1_to_M,N} and \ref{from_k=1_to_kinN} to the result that the analogon of the difference in Theorem \ref{cubatureBermudanconv} is of order $s^2\cdot 2^{-M}$, too:

\begin{Th} \label{cubatureBermudanconv_barf=min} Suppose $\bar f=\max_{j\in\{1,\dots,d\}}\exp\left((\cdot)_j\right)$ and consider a compact set $B$ such that $B-(\ln K)_{j=1}^d\subseteq\left[-{R},R\right]^d$. Assume that $y_{i}^{(h)}\geq 0$ for all $i\in\{1,\dots,m\}$ and $0\in\left\{y_i^{(h)} \ : \ i\in\{1,\dots,m\}\right\}$. Under these assumptions there exists a real number $D>0$, given by $$D=\left(\left(\ln\gamma_1-r\right)\tilde{D} +rK+ C_0  \right)\cdot R^{d-1}\cdot \frac{\sum_{j=1}^d\max_{i} \left(\left(y_i^{(h)}\right)_j\vee 0\right)}{h} $$ (with $C_0$ and $\tilde D$ as defined in Lemma \ref{estimateonE} and the $\gamma_1$ of Lemma \ref{barf=max_gamma1}) such that for all $N>M\in\NN$, $s\in(0,T)\cap \left(2^N\cdot I\right)$ and $f\geq g\vee 0$, one has 
\begin{eqnarray*}&&\left\|\left(B_{s\cdot 2^{-N}}\right)^{\circ\left(2^{N}\right)}f - \left(B_{s\cdot 2^{-M}}\right)^{\circ\left(2^M\right)}f\right\|_{L^1\left(E^h \cap \bigcap_{\ell}\left(B-y_\ell^{\left(s\right)}\right)\right)} \\ &\leq& D\cdot {s}^2\cdot{2}^{-M}\left(1- 2^{-(N-M-1)}\right) \\ &\leq& D\cdot {s}^2\cdot{2}^{-M}\longrightarrow 0 \text{ as }M\rightarrow\infty. \end{eqnarray*}
\end{Th}

\begin{rem} Suppose one is looking at a market model where $\min_{i\in\{1,\dots,m\}}y_{i}^{(h)}$ is still positive. Then one has \begin{align}\label{B^Qbound}\forall x\in\RR^d\quad \left(B_{s\cdot 2^{-N}}\right)^{\circ\left(2^{N}\right)}(g\vee 0)(x)\leq \left(B^{Q}_{s\cdot 2^{-N}}\right)^{\circ\left(2^{N}\right)}(g\vee 0)(x)\end{align} when taking $Q$ to be the translation-invariant Markov chain with increments $y^{(h)}_1-\min_{i\in\{1,\dots,m\}}y_{i}^{(h)},\dots, y^{(h)}_m-\min_{i\in\{1,\dots,m\}}y_{i}^{(h)}$ and weights $\alpha^{(h)}_1,\dots, \alpha^{(h)}_m $ (that is the Markov-chain obtained from $P$ by simply shifting the increments in the componentwise negative direction $-\min_{i\in\{1,\dots,m\}}y_{i}^{(h)}$) respectively, and defining the operator family $\left(B^Q\right)_{t\in I}$ from the semigroup $Q$ via $B_t^Q:f\mapsto e^{-rt}Q_tf\vee g$ for all $t\in I$ (viz. just as the family $B$ was defined from the semigroup $P$). Note that just as in the proof of Remark \ref{treeispolynomial} one can use deduce upper bounds on the increments of $P_nh$ from the regularity of the increments of $P_h$ (denoted by $\xi_1,\dots,\xi_m$ in Remark \ref{treeispolynomial}; and this reasoning is applicable to the Markov chain $Q$ whenever it is applicable to $P$. Also note that the upper bound in estimate \eqref{B^Qbound} is sharp in the sense that the right hand side approaches the value of the left hand side as $\min_{i\in\{1,\dots,m\}}y_{i}^{(h)}$ approaches zero. 
\end{rem}

\begin{rem} This $L^1$-convergence result has some practical interest, as in practice quite frequently the exact start price of the (multiple) asset on which an option is issued, is unknown. Instead, one will have the logarithmic start price vector $x\in\RR^d$ a short time $\delta>0$ before the actual option contract becomes valid. Now, assuming that $\PP^x_{X_\delta}$ has a continuous density $\frac{\PP^x_{X_\delta}}{\lambda^d}$, this function $\frac{\PP^x_{X_\delta}}{\lambda^d}$ will be bounded on $E^h \cap \bigcap_{\ell}\left(B-y_\ell^{\left(s\right)}\right)$ by some constant $$C:=\sup_{E^h \cap \bigcap_{\ell}\left(B-y_\ell^{\left(s\right)}\right)} \frac{\PP^x_{X_\delta}}{\lambda^d}<+\infty.$$ One will therefore have for all $f\geq g\vee 0$, $s\in(0,T]\cap I$ and $N>M\in\NN$, \begin{eqnarray*}&&\EE^x\left[\begin{array}{c}\left|\left(B_{s\cdot 2^{-N}}\right)^{\circ\left(2^{N}\right)}f - \left(B_{s\cdot 2^{-M}}\right)^{\circ\left(2^M\right)}f\right|\left(X_\delta\right),\\ X_\delta\in E^h \cap \bigcap_{\ell}\left(B-y_\ell^{\left(s\right)}\right)\end{array}\right]\\&\leq&
\left\|\left(B_{s\cdot 2^{-N}}\right)^{\circ\left(2^{N}\right)}f - \left(B_{s\cdot 2^{-M}}\right)^{\circ\left(2^M\right)}f\right\|_{L^1\left(\PP^x_{X_\delta}\left[E^h \cap \bigcap_{\ell}\left(B-y_\ell^{\left(s\right)}\right)\cap\cdot\right]\right)} \\ &\leq& C\cdot\left\|\left(B_{s\cdot 2^{-N}}\right)^{\circ\left(2^{N}\right)}f - \left(B_{s\cdot 2^{-M}}\right)^{\circ\left(2^M\right)}f\right\|_{L^1\left(\lambda^d\left[E^h \cap \bigcap_{\ell}\left(B-y_\ell^{\left(s\right)}\right)\cap\cdot\right]\right)} \\ &\leq& C\cdot D\cdot {s}^2\cdot{2}^{-M}. \end{eqnarray*}

\end{rem}

{\bf Acknowledgements.} The author would like to thank the German Academic Exchange Service for the pre-doctoral research grant he received ({\em Doktorandenstipendium des Deutschen Akademischen Austauschdienstes}) and the German National Academic Foundation ({\em Studienstiftung des deutschen Volkes}) for their generous support in both financial and non-material terms. 

Moreover, he owes a huge debt of gratitude to his supervisor, Professor Terry J Lyons, for numerous extremely helpful discussions.

\end{document}